\documentclass{elsart}

\usepackage{chemarr}
\usepackage{amsfonts}
\usepackage{amssymb}
\usepackage{latexsym}
\usepackage{graphicx}
\usepackage{hyperref}
\usepackage{amsbsy}
\usepackage{psfrag}
\usepackage{color}

\renewcommand{\And}{\mbox{\ and \ }} 

\newcommand{\D}{\displaystyle}
\newcommand{\restr}[1]{\raisebox{-0.3em}{$\lb|_{#1}\rb.$}} 
\newcommand{\ignore}[1]{}    
\newcommand{\breath}{\medskip} 
\newtheorem{THM}{Theorem}[section]

\newcounter{claimcount}[THM]  
\newcounter{subclaimcount}[claimcount]
\newtheorem{PROP}[THM]{Proposition} 

\newtheorem{lemma}[THM]{Lemma} 

\newtheorem{COR}[THM]{Corollary}

\newcommand{\dfn}{\sf\em} 
\newcommand{\Theorem}[2]{\begin{THM}{\sf #1}  #2 \end{THM}}
\newcommand{\Proposition}[2]{\begin{PROP}{\sf #1}  #2 \end{PROP}}
\newcommand{\Lemma}[2]{\begin{lemma}{\sf #1}  #2 \end{lemma}}
\newcommand{\Corollary}[2]{\begin{COR}{\sf #1}  #2 \end{COR}} 
\newcommand{\THMfont}[1]{{\sl #1}}      
\newcommand{\Definition}[2]{ \refstepcounter{THM} {\bf Definition \theTHM:} {\em #1} #2  \hfill $\diamondsuit$}  

\newcommand{\example}[1]{    
    \refstepcounter{THM}    
                 \begin{list}{} 
			{\setlength{\leftmargin}{0em} 
			\setlength{\rightmargin}{0em}}   
     \item {\bf Example \theTHM:} #1                 
  \hfill$\diamondsuit$  \end{list}   			

}     
\newcommand{\bthmlist}{ \begin{list}{{\bf(\alph{enumi})}} 
{\usecounter{enumi}
 \setlength{\itemsep}{0.2em}
 \setlength{\topsep}{0.2em}
 \setlength{\itemindent}{0em} 
\setlength{\parsep}{0em} 
} 
} 

\newcommand{\ethmlist}{\end{list}}    
\newcommand{\Claim}[1]{\refstepcounter{claimcount} \vspace{0.3em}               \noindent {\sc Claim \theclaimcount: \ }\THMfont{ #1}} 
\newcommand{\Subclaim}[1]{\refstepcounter{subclaimcount} \vspace{0.3em}               \noindent {\sc Claim \theclaimcount.\thesubclaimcount: \ }\THMfont{ #1}} 
\newcommand{\subclaim}{\Subclaim}

\newcommand{\bprf}[1][Proof:]{\begin{list}{}
 			{\setlength{\leftmargin}{0.5em}
 			\setlength{\rightmargin}{0em}
 			\setlength{\listparindent}{1em}}
                         \item {\em \hspace{-1em}  #1  }} 

\newcommand{\eprf}{\end{list}} 
\newcommand{\bthmprf}{\bprf}
\newcommand{\bclaimprf}{\bprf}
\newcommand{\bsubclaimprf}{\bprf}
\newcommand{\ethmprf}{ \hfill$\Box$  \eprf  \breath  } 
\newcommand{\eclaimprf}{ \hfill $\Diamond$~{\scriptsize {\tt Claim~\theclaimcount}}\eprf}  
\newcommand{\esubclaimprf}{ \hfill $\triangledown$~{\scriptsize  {\tt Claim~\theclaimcount.\thesubclaimcount}}\eprf}  
\newcommand{\beq}{\begin{eqnarray*}}
\newcommand{\eeq}{\end{eqnarray*}} 
\newcommand{\beqn}{ \begin{equation} }
\newcommand{\eeqn}{ \end{equation} }
\newcommand{\blist}{\begin{enumerate}}
\newcommand{\elist}{\end{enumerate}} 
\newcommand{\bitem}{\begin{itemize}}
\newcommand{\eitem}{\end{itemize}} 
\newcommand{\If}{\mbox{\ if \ }}  
\newcommand{\done}{{\mathsf{ 1\!\!1}}} 
\newcommand{\dI}{{\mathbb{I}}}
\newcommand{\dJ}{{\mathbb{J}}}
\newcommand{\dN}{{\mathbb{N}}}
\newcommand{\dU}{{\mathbb{U}}}
\newcommand{\dV}{{\mathbb{V}}}
\newcommand{\dW}{{\mathbb{W}}}
\newcommand{\dX}{{\mathbb{X}}}
\newcommand{\dY}{{\mathbb{Y}}}
\newcommand{\dZ}{{\mathbb{Z}}}       %  Define Bar macros 
\newcommand{\barmu}{{\overline{\mu }}}
\newcommand{\barba}{{\overline{\mathbf{ a}}}}
\newcommand{\bA}{{\mathbf{ A}}}
\newcommand{\bB}{{\mathbf{ B}}}
\newcommand{\bC}{{\mathbf{ C}}}
\newcommand{\bD}{{\mathbf{ D}}}
\newcommand{\bE}{{\mathbf{ E}}}
\newcommand{\bG}{{\mathbf{ G}}}
\newcommand{\bL}{{\mathbf{ L}}}
\newcommand{\bR}{{\mathbf{ R}}}
\newcommand{\bS}{{\mathbf{ S}}}
\newcommand{\bU}{{\mathbf{ U}}}
\newcommand{\bX}{{\mathbf{ X}}}
\newcommand{\ba}{{\mathbf{ a}}}
\newcommand{\bb}{{\mathbf{ b}}}
\newcommand{\bc}{{\mathbf{ c}}}
\newcommand{\bd}{{\mathbf{ d}}}
\newcommand{\bl}{{\mathbf{ l}}} 
\newcommand{\br}{{\mathbf{ r}}}
\newcommand{\bs}{{\mathbf{ s}}}
\newcommand{\bt}{{\mathbf{ t}}}
\newcommand{\bx}{{\mathbf{ x}}}
\newcommand{\by}{{\mathbf{ y}}}
\newcommand{\bz}{{\mathbf{ z}}}    \usepackage{amsbsy} 
\newcommand{\sA}{{\mathcal{ A}}}
\newcommand{\sB}{{\mathcal{ B}}}
\newcommand{\sC}{{\mathcal{ C}}}
\newcommand{\sD}{{\mathcal{ D}}}
\newcommand{\sF}{{\mathcal{ F}}}
\newcommand{\sI}{{\mathcal{ I}}}
\newcommand{\sL}{{\mathcal{ L}}}
\newcommand{\sM}{{\mathcal{ M}}}
\newcommand{\sO}{{\mathcal{ O}}}
\newcommand{\sP}{{\mathcal{ P}}}
\newcommand{\sR}{{\mathcal{ R}}}
\newcommand{\sT}{{\mathcal{ T}}}
\newcommand{\sX}{{\mathcal{ X}}}
\newcommand{\alp }{\alpha}
\newcommand{\bet }{\beta}
\newcommand{\gam }{\gamma}
\newcommand{\del }{\delta}
\newcommand{\lam }{\lambda}
\newcommand{\sig }{\sigma} 
\newcommand{\omg }{\omega}
\newcommand{\Gam }{\Gamma}
\newcommand{\Del }{\Delta}
\newcommand{\Ups }{\Upsilon}
\newcommand{\Omg }{\Omega}   %  Define Hat macros 
\newcommand{\h}{\widehat} 
\newcommand{\hz}{{\widehat{z}}}   
\newcommand{\hbL}{{\widehat{\mathbf{ L}}}}
\newcommand{\hbR}{{\widehat{\mathbf{ R}}}}
\newcommand{\hba}{{\widehat{\mathbf{ a}}}}
\newcommand{\hsA}{{\widehat{\mathcal{ A}}}}
\newcommand{\hlam }{{\widehat{\lambda}}}
\newcommand{\hUps }{{\widehat{\Upsilon}}}
\newcommand{\hPhi}{{\widehat{\Phi }}}
\newcommand{\fA}{{\mathsf{ A}}}
\newcommand{\fD}{{\mathsf{ D}}}
\newcommand{\fF}{{\mathsf{ F}}}
\newcommand{\fM}{{\mathsf{ M}}}
\newcommand{\fP}{{\mathsf{ P}}}
\newcommand{\fT}{{\mathsf{ T}}}
\newcommand{\fv}{{\mathsf{ v}}}
\newcommand{\tl}{\widetilde} 
\newcommand{\tlbA}{{\widetilde{\mathbf{ A}}}}
\newcommand{\tlbS}{{\widetilde{\mathbf{ S}}}}
\newcommand{\vV}{{\vec{V}}}
\newcommand{\vv}{{\vec{v}}}
\newcommand{\lb}{\left}
\newcommand{\rb}{\right} 
\newcommand{\Array}[2][cccccccccccccccccccccccccccccccccccc] {{\begin{array}{#1}#2\end{array}}}      
\newcommand{\maketall}{\rule[-0.5em]{0em}{1em}}        
\newcommand{\map}{{\longrightarrow}}
\newcommand{\goto}{{\rightarrow}}
\newcommand{\into}{{\map}}
\newcommand{\statement}[1]{\lb(  \maketall       \begin{minipage}{40em}       \begin{tabbing}         #1        \end{tabbing}      \end{minipage}  \rb)}     
\newcommand{\oo}{{\infty}}        
\newcommand{\X}{\times}
\newcommand{\x}{\X}
\newcommand{\tensor}{\otimes}
\newcommand{\union}{\cup}
\newcommand{\Union}{\bigcup}
\newcommand{\intsct}{\cap}
\newcommand{\disj}{\sqcup}
\newcommand{\Disj}{\bigsqcup}   
\newcommand{\set}[2]{{\left\{ #1 \; ; \; #2 \right\} }} 
\newcommand{\Id}[1]{{\mathbf{ Id}_{{#1}}}}
\newcommand{\pr}[1]{{\mathbf{ pr}_{{#1}}}}
\newcommand{\chr}[1]{{{\done}_{{#1}}}} 
\newcommand{\choice}[1]{{\lb\{ \begin{array}{rcl}                                 #1                                \end{array}  \rb.  }}                     %        
\newcommand{\eeequals}[1]{\raisebox{-0.9ex}{$\overline{\overline{{\scriptscriptstyle{\mathrm{#1}}}}}$}} 
\newcommand{\Fix}[1]{{\sf Fix}\lb[#1\rb]}
\newcommand{\shift}[1]{\sig^{#1}}    
\newcommand{\End}[2][]{\mathsf{End}_{#1} \lb(#2\rb)}
\newcommand{\Natur}{\dN}
\newcommand{\Zahl}{\dZ}
\newcommand{\Zahlmod}[1]{{\Zahl_{/#1}}}
\newcommand{\CC}[1]{{\lb[ #1 \rb]}}
\newcommand{\CO}[1]{{\lb[ #1 \rb)}}
\newcommand{\OC}[1]{{\lb( #1 \rb]}}
\newcommand{\OO}[1]{{\lb( #1 \rb)}}   
\newcommand{\gB}{\bB}
\newcommand{\gE}{\bE}
\newcommand{\gG}{\bG}
\newcommand{\gL}{\bL}
\newcommand{\gR}{\bR}
\newcommand{\ZD}[1][D]{{\Zahl^{#1}}}
\newcommand{\AZD}[1][D]{\sA^{\ZD[#1]}}
\newcommand{\AZ}{\sA^{\Zahl}}
\newcommand{\AN}{\sA^{\Natur}}
\newcommand{\BN}{\sB^{\Natur}}
\newcommand{\BZ}{\sB^{\Zahl}}
\newcommand{\RZ}{\sR^{\Zahl}}
\newcommand{\LZ}{\sL^{\Zahl}}
\newcommand{\XN}{\sX^{\Natur}}
\newcommand{\black}{{\scriptstyle{\blacksquare}}}
\newcommand{\white}{{\scriptstyle{\square}}}  
\newcommand{\NRemark}[1]{    
                 \begin{list}{}
 			{\setlength{\leftmargin}{0em}
 			\setlength{\rightmargin}{0em}}
        \item        \refstepcounter{THM} {\em Remark \theTHM:} #1
                   \hfill$\diamondsuit$  \end{list}   		

	} 
\newcommand{\Remarks}[1]{
                     \begin{list}{}
 			{\setlength{\leftmargin}{0em}
 			\setlength{\rightmargin}{0em}}
        \item {\em Remarks:} #1
                   \hfill$\diamondsuit$  \end{list}

   			} 
\newcommand{\Remark}[1]{
                     \begin{list}{}
 			{\setlength{\leftmargin}{0em}
 			\setlength{\rightmargin}{0em}}
        \item {\em Remark:} #1
                   \hfill$\diamondsuit$  \end{list}
   			}
 
\newcommand{\hgL}{\widehat{\gL}}
\newcommand{\hgR}{\widehat{\gR}} 
\newcommand{\Turing}{\fT\!\fM}
\newcommand{\APDA}{\fA\!\fP\!\fD\!\fA}
\newcommand{\AFA}{\fA\!\fF\!\fA}  
\newcommand{\ECA}[1]{{}_{\varepsilon}\!\Phi_{\mbox{\tiny #1}}}
\newcommand{\lcm}{\mathrm{lcm}}    
\newcommand{\muae}{{\rm ($\mu$-\ae)}}   
\newcommand{\SigAlg}{{\mathfrak{S}}}  
\newcommand{\parag}[1]{\vspace{0.3em} {\em #1}}  
\newcommand{\backtau}{\stackrel{\leftarrow}{\tau}}   \DeclareMathAlphabet{\zapf}{OT1}{pzc}{m}{it}  
\newcommand{\za}{{\zapf{a}}}
\newcommand{\zb}{{\zapf{b}}}
\newcommand{\zc}{{\zapf{c}}}
\newcommand{\zd}{{\zapf{d}}}
\newcommand{\ze}{{\zapf{e}}}
\newcommand{\zi}{{\zapf{i}}}
\newcommand{\zl}{{\zapf{l}}}
\newcommand{\zo}{{\zapf{o}}}
\newcommand{\zp}{{\zapf{p}}}
\newcommand{\zr}{{\zapf{r}}}
\newcommand{\zt}{{\zapf{t}}}
\newcommand{\zx}{{\zapf{x}}}
\newcommand{\zy}{{\zapf{y}}}
\newcommand{\hza}{{\widehat{\zapf{a}}}}
\newcommand{\hzd}{{\widehat{\zapf{d}}}}
\newcommand{\hzl}{{\widehat{\zapf{l}}}}
\newcommand{\hzr}{{\widehat{\zapf{r}}}}
\newcommand{\TM}{\mbox{\large{\tt M}}}   
\newcommand{\Vdrift}{\vV_{\scriptscriptstyle{\mathrm{drift}}}}  
\renewcommand{\gB}{\bB}
\renewcommand{\gE}{\bE}
\renewcommand{\gG}{\bG}
\renewcommand{\gL}{\bL}
\renewcommand{\gR}{\bR}

\begin{document}

\begin{frontmatter}
\title{Defect Particle Kinematics in One-Dimensional Cellular Automata}
\author{Marcus Pivato\corauthref{myaddress}}
\address{Dept. of Mathematics \& Computer Science, Wesleyan University \\ and Department of Mathematics,  Trent University}

\ead{pivato@xaravve.trentu.ca, marcuspivato@trentu.ca}
\corauth[myaddress]{Trent University, 1600 West Bank Drive, Peterborough, Ontario, Canada, K9J 7B8}

%spellchecked
\begin{abstract}
Let $\AZ$ be the Cantor space of bi-infinite sequences in a finite
alphabet $\sA$, and let $\shift{}$ be the shift map on
$\AZ$.  A {\dfn cellular automaton} is a continuous,
$\shift{}$-commuting self-map $\Phi$ of $\AZ$, and a {\dfn
$\Phi$-invariant subshift} is a closed, $(\Phi,\shift{})$-invariant
subset $\bS\subset\AZ$.  Suppose $\ba\in\AZ$ is $\bS$-admissible
everywhere except for some small region we call a {\dfn defect}. It
has been empirically observed that such defects persist under
iteration of $\Phi$, and often propagate like `particles'.  We
characterize the motion of these particles, and show that it falls
into several regimes, ranging from simple deterministic motion, to
generalized random walks, to complex motion emulating Turing machines
or pushdown automata.  One consequence is that some questions
about defect behaviour are formally undecidable.

\end{abstract}


\end{frontmatter}

  A recurring theme in cellular automata is the emergence and
persistence of homogeneous `domains' (each characterized by a
particular spatial pattern), separated by {\em defects} (analogous to
`domain boundaries' or `kinks' in a crystalline solid).
 Defects were first empirically observed by Grassberger in the
`elementary' cellular automata or `ECA' (radius-one CA on
$\{0,1\}^\Zahl$) with numbers \#18,
\#122, \#126, \#146, and \#182 \cite{Gra83,Gra84} and also noted in ECA \#184,
which was originally studied as a simple model of surface growth
\cite[\S III.B]{KrSp88}, and later as a model of single-lane traffic  
\cite{Fuks99,Blank03,BelitskyFerrari05}.
Based on Grassberger's
observations, Lind \cite[\S5]{Lin84} conjectured that the defects of
ECA\#18 perform random walks.  This conjecture was reiterated by
Boccara {\em et al.}, who empirically investigated the motion and
interactions of defects in ECA
\#18 and also \#54, \#62, and \#184 (see Figure \ref{fig:defect.intro}),
and longer range totalistic CA \cite{BoRo91,BNR91}; see also 
\cite[\S3.1.2.2 \& \S3.1.4.4]{Ilachinski}.

Eloranta developed the first rigorous mathematical theory of cellular
automaton defects in
\cite{Elo93a,Elo93b,Elo94,Elo95}, and, together with
Numelin, proved Lind's conjecture in \cite{ElNu}.  Meanwhile,
Crutchfield and Hanson developed an empirical methodology called {\em
Computational Mechanics}
\cite{Han}, which they  applied to ECA\#18
\cite{CrHa92,CrHa93a} and other CA contrived to act like
ECA\#18 \cite{CrHa93b}, as well as ECA\#54 \cite{CrHa97,CrHR} and
ECA\#110 \cite{CrHR}.  They  also obtained a
tight theoretical bound on the number of possible reactions between
two defects \cite{CrHR} (improving an earlier result of
\cite{PaStTh}).  Finally, using genetic algorithms, they and their
collaborators `bred' CA which performed computations such as
synchronization or density-classification.  A careful analysis then
revealed that these CA performed their computations through
propagating and interacting defects;   this `particle-based
computation' had emerged spontaneously through natural selection
\cite{CDM,CrHM,CrHR}.

\begin{figure}[h]
\centerline{\begin{tabular}{l}
\includegraphics[clip=true, trim=0 0 350 30, width=33em,height=13em]{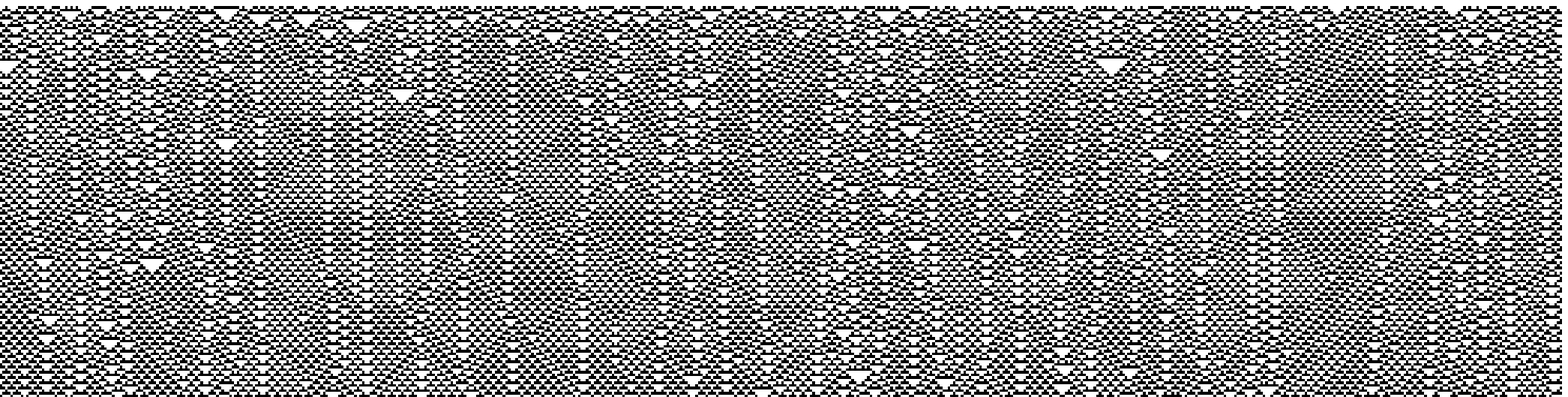} \\
 {\footnotesize{\bf(A)} \ ECA \#54 } \\
\includegraphics[clip=true, trim=0 30 350 0, width=33em,height=13em]{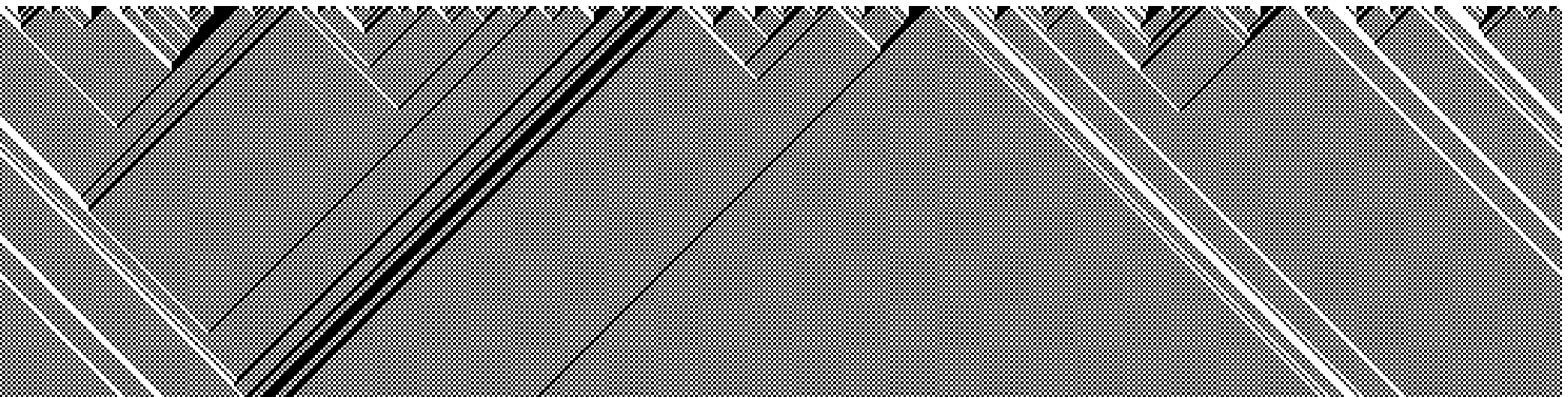} \\
{\footnotesize{\bf(B)} \ ECA \#184} \\
\end{tabular}}
\caption{{\footnotesize
 Spacetime diagrams showing defect dynamics in two
cellular automata.  Each picture show 120 timesteps on a 300 pixel array
(time increases downwards).
 \label{fig:defect.intro}}}
\end{figure}

\ignore{
\begin{figure}[h]
\centerline{\begin{tabular}{ll}
\includegraphics[clip=true, trim=0 0 350 30, width=16.5em,height=6.5em]{Pictures/rule54a.eps} 
& \includegraphics[clip=true, trim=0 30 350 0, width=16.5em,height=6.5em]{Pictures/rule184a.eps} \\
 {\footnotesize{\bf(A)} \ ECA \#54 } & {\footnotesize{\bf(B)} \ ECA \#184} \\
\end{tabular}}
\caption{{\footnotesize
 Spacetime diagrams showing defect dynamics in two
cellular automata.  Each picture show 120 timesteps on a 300 pixel array
(time increases downwards).
 \label{fig:defect.intro}}}
\end{figure}
}
  In two companion papers \cite{PivatoDefect2,PivatoDefect1}, we
develop algebraic invariants to explain why defects persist under
iteration of the cellular automaton, instead of disappearing.  These
defects often behave like `particles', which propagate through space
until they collide and interact with other defects.  In this paper, we
characterize the motion of these `defect particles', when the
background domain is a one-dimensional subshift of finite type which
is invariant under the action of a one-dimensional cellular automaton.
In \S\ref{S:particle.prop} we formally define `defect particles' and
introduce a framework to investigate their motion.  Depending on the
$(\Phi,\shift{}$)-dynamical properties of the ambient subshift, the
defect particle falls into one of several `kinematic regimes', ranging
from ballistic motion (\S\ref{S:ballistic}), to a generalized random
walk (\S\ref{S:diffuse}), to the emulation of Turing machines or
pushdown automata (\S\ref{S:turing}).  Sections \S\ref{S:ballistic},
\S\ref{S:diffuse} and \S\ref{S:turing} are logically independent of
one another.

\subsubsection*{Preliminaries \& Notation}

For any $L\leq R\in\Zahl$, we define $\CC{L...R}:=\{L, L\!+\!1,\,\ldots,R\}$,
\ $\CO{L...R}:=\CC{L\ldots R\!-\!1}$, \ $\OC{L...R}:=\CC{L\!+\!1\ldots
R}$, etc.  We likewise define $\OC{-\oo....R}$, \ $\CO{L...\oo}$, etc.
  Let $\sA$ be a finite alphabet, and
let $\AZ$ be the set of all doubly-infinite sequences in $\sA$,
which we write as
$\ba=[\za_{ z}]_{z\in\Zahl}$, where $\za_z\in\sA$ for all $z\in\Zahl$.
The {\dfn Cantor metric} on $\AZ$ is defined by
$d(\ba,\bb)=2^{-\Del(\ba,\bb)}$, where
$\Del(\ba,\bb):=\min\set{| z|}{\za_z\neq \zb_z}$.  It follows that
$(\AZ,d)$ is a Cantor space (i.e. a compact, totally disconnected,
perfect metric space).  If $\ba\in\AZ$, and $\dU\subset\Zahl$,
then we define $\ba_\dU\in\sA^\dU$ by
$\ba_\dU:=[\za_u]_{u\in\dU}$.  If $ z\in\Zahl$, then strictly speaking,
$\ba_{z+\dU}\in\sA^{ z+\dU}$; however, it is sometimes convenient
to `abuse notation' and treat $\ba_{ z+\dU}$ as an element of
$\sA^{\dU}$ in the obvious way.

\ignore{  If $\dX\subset\dY\subseteq\Zahl$, and
$\bx\in\sA^\dX$ and $\by\in\sA^\dY$, we say $\bx$ {\dfn occurs} in
$\by$ (``$\bx \sqsubset \by$'') if $\bx=\by_{\dX}$.}

\ignore{  If $\ba=\za_1\za_2\ldots \za_p\in\sA^\CC{1..p}$, then we define the
$p$-periodic sequence $\barba:=[\ldots  \za_1\za_2\ldots
\za_p\ \underscore{\za_1}\za_2\ldots \za_p \ \za_1\za_2\ldots \za_p\ldots]$ (we underline the
zeroth entry).  More generally, if $\dU\subset\Zahl$ is any subset which
`tiles' $\Zahl$ in an obvious fashion (e.g. a square), and
$\ba\subset\sA^\dU$, then we define the periodic extension 
$\barba\in\AZ$ of $\ba$ in the obvious sense.}

We define the {\dfn shift map} $\shift{}:\AZ\into\AZ$ by
$\shift{}(\ba)_{z} = \za_{z+1}$ for all $\ba\in\AZ$ and
$ z\in\Zahl$.  A {\dfn cellular automaton} is a
transformation $\Phi:\AZ\into\AZ$ that is
continuous and commutes with $\shift{}$.
Equivalently, $\Phi$ is determined by
a {\dfn local rule} $\phi:\sA^\CC{-r...r}\into\sA$ (for some $r\in\Natur$)
such that
$\Phi(\ba)_{ z} = \phi(\ba_{\CC{z\!-\!r\ldots z\!+\!r}})$ for all $\ba\in\AZ$ and
$z\in\Zahl$ \cite{Hedlund};  we  say that $\Phi$ has {\dfn radius} $r$.

  A subset $\bS\subset\AZ$ is a {\dfn subshift}
\cite{LindMarcus,Kitchens} if $\bS$ is closed in the Cantor topology,
and $\shift{}(\bS)=\bS$.  For any
$\dU\subset\Zahl$, we define $\bS_\dU:=\set{\bs_\dU}{\bs\in\bS}$. 
In particular, for any $q>0$, let 
$\bS_q := \bS_{\CO{0...q}}$ be the set of {\dfn admissible
$q$-words} for $\bS$.  We say $\bS$ is {\dfn subshift of
finite type} (SFT) if there is some $q>0$ (the {\dfn radius} of $\bS$)
such that  $\bS$ is entirely described by $\bS_q$, in the sense
that $\bS=\set{\bs\in\AZ}{\bs_{\CO{z...z\!+\!q}}\in\bS_q, \
\forall z\in\Zahl}$.  If $q=2$, then $\bS$ is called a
{\dfn Markov subshift}, and the elements of
$\bS_{2}\subseteq\sA^2$ are called {\dfn admissible transitions}; 
equivalently, $\bS$ is the set of all bi-infinite directed paths in a 
digraph  whose vertices are the elements of $\sA$, with
an edge $\za\leadsto \zb$ iff $(\za,\zb)\in\bS_2$. 

  If $\Phi:\AZ\into\AZ$ is a cellular automaton, then we say $\bS$
is (weakly) {\dfn $\Phi$-invariant} if $\Phi(\bS)\subseteq\bS$
(i.e. $\Phi$ is an {\dfn endomorphism} of $\bS$).  For example,
the set $\Fix{\Phi}:=\set{\ba\in\AZ}{\Phi(\ba)=\ba}$  of
$\Phi$-{\dfn fixed points} is a $\Phi$-invariant
SFT.  Likewise, if $p\in\Natur$ and $v\in\Zahl$, then the set
$\Fix{\Phi^p}$ of  {\dfn $(\Phi,p)$-periodic points} and the set
$\Fix{\Phi^p\circ\shift{-pv}}$ of {\dfn $(\Phi,p,v)$-travelling waves} are
$\Phi$-invariant SFTs.  Also, for any $p\in\Natur$, the set
$\Fix{\shift{p}}$ of $p$-periodic sequences is a $\Phi$-invariant
SFT.

 If $\Phi$ has radius $r$, then for any $q>0$,
$\Phi$ induces a function $\Phi:\sA^{q+2r}\into\sA^{q}$.
If $\bS\subset\AZ$ is an SFT determined by a set $\bS_q\subset\sA^{q}$ of
admissible $q$-blocks, then $\statement{$\Phi(\bS)\subseteq\bS)$}
\iff \statement{$\Phi(\bS_{q+2r}) \subseteq \bS_q$}$. 
The monoid of endomorphisms of an SFT can be quite huge;
see \cite[Ch.3]{Kitchens} or \cite[\S13.2]{LindMarcus}.

If $\bS\subset\AZ$ is a subshift, then we define $\bS^- :=\bS_\OC{-\oo...-1}
\subseteq \sA^\OC{-\oo...-1}$  to be the set of all left-infinite,
$\bS$-admissible sequences, and define
$\bS^+ :=\bS_\CO{1...\oo}
\subseteq \sA^\CO{1...\oo}$  to be the set of all right-infinite,
$\bS$-admissible sequences.

\breath

{\em Notation \& Font conventions:} Upper case calligraphic letters
($\sA,\sB,\sC,\ldots$) denote finite symbolic alphabets (of cellular
automata, Turing machines, etc.).  Upper-case boldface
letters ($\bA,\bB,\bC,\ldots$) denote subsets of $\AZ$
(e.g. subshifts);  lowercase bold-faced letters ($\ba,\bb,\bc,\ldots$)
denote elements of $\AZ$. Zapf letters ($\za,\zb,\zc,\ldots$) are
elements of $\sA$;  Roman letters ($a,b,c,\ldots$) are integers.
Upper-case hollow font ($\dU,\dV,\dW,\ldots$) are subsets of $\Zahl$,
upper-case Greek letters ($\Phi,\Psi,\ldots$) denote functions on
$\AZ$ (e.g. CA), and lower-case Greek letters ($\phi,\psi,\ldots$)
denote other functions (e.g. local rules, probability measures).

  We generally indicate related objects by related letters.
For example, if $\sL,\sR\subset\sA$ are two subalphabets, then a
subshift of $\LZ$ would be denoted by $\bL$, with typical element
$\bl:=[\zl_z]_{z\in\Zahl}\in\bL$ (where $\zl_z\in\sL$), whereas 
a subshift of $\RZ$ would be denoted by $\bR$, with typical element
$\br:=[\zr_z]_{z\in\Zahl}\in\bR$ (where $\zr_z\in\sR$).

\section{Defect Particles\label{S:particle.prop}}

  Let $\bS\subset\AZ$ be a subshift of finite type, and let
$\Phi:\AZ\into\AZ$ be a one-dimensional cellular automaton with
$\Phi(\bS)\subseteq\bS$.  By passing to a `higher block presentation',
we can assume that $\Phi$ is a nearest-neighbour CA and that $\bS$ is
a Markov subshift.  To be precise, suppose $\Phi$ has radius $r$ and
that $\bS$ is determined by a set $\bS_{q}$ of admissible $q$-blocks.
Let $P:=\max\{2r,q\}$, let $\sB:=\sA^P$, and let $\tlbA\subset\BZ$ be
the {\dfn $P$-block presentation} of $\AZ$ (see
\cite[Defn.1.4.1]{LindMarcus} or \cite[Fig.1.4.1]{Kitchens}).  That
is, $\tlbA$ is the the Markov subshift defined by the digraph with
vertex set $\sA^P$, with an edge
$[\za_1,\ldots,\za_P]\leadsto[\zb_1,\ldots,\zb_P]$ iff $\zb_p=\za_{p+1}$ for all
$p\in\CO{1...P}$ (this is sometimes called the {\dfn de~Bruijn digraph} of
$\sA^P$).  Thus, $\Phi$ is conjugate to an endomorphism
$\tl\Phi:\tlbA\into\tlbA$, which can be extended (in an arbitrary way) to
a cellular automaton $\Psi:\BZ\into\BZ$ such that
$\Psi(\tlbA)\subseteq\tlbA$ and $\Psi\restr{\tlbA}=\tl\Phi$.  Let $\tlbS$ be
the image of $\bS$ inside $\tlbA$; then $\tlbS$ is a Markov subshift of
$\BZ$, and $\Psi(\tlbS)\subseteq\tlbS$.  Now replace $\bS$ with $\tlbS$ and $\Phi$
with $\Psi$.

   Let $\ba^0\in\AZ$, then $\ba^0$ has a {\dfn
single defect} in the interval $\CC{i...k}\subset\Zahl$ if
$[\za_j,\za_{j+1}]\in\bS_{2}$ for all
$j\not\in\CC{i...k}$, while $[\za_j,\za_{j+1}]\not\in\bS_{2}$ for all
$j\in\CC{i...k}$.  If $i=k$, then the defect has {\dfn width $0$}, and
consists of a single inadmissible transition between two
half-infinite, $\bS$-admissible sequences: \beqn
\label{zero.width.defect}
\ba^0\quad=\quad
 [\underbrace{\underline{\ldots \za^0_{i-3} , \za^0_{i-2} , \za^0_{i-1} \za^0_i}}_{\scriptscriptstyle\mathrm{admissible}} \!\!\raisebox{-1.4em}{${\Uparrow}\atop{\scriptscriptstyle\mathrm{defect}}$}\!\!  \underbrace{\underline{\za^0_{i+1} , \za^0_{i+2} , \za^0_{i+3}\ldots}}_{\scriptscriptstyle\mathrm{admissible}}]
\eeqn
(here we underline the admissible sequences for visibility).
If $i<k$, then we say that $(\za_{i+1},\ldots,\za_{k})$ is a {\dfn defect word} of
{\dfn width} $w:=k-i$:
\beqn
\label{variable.width.defect.0}
\ba^0\quad=\quad
 [\underline{\ldots \za^0_{i-3} , \za^0_{i-2} , \za^0_{i-1} \za^0_i} \quad 
\underbrace{ \za^0_{i+1} \ldots \za^0_{k}}_{\scriptscriptstyle\mathrm{defect}}
 \quad 
 \underline{\za^0_{k+1} , \za^0_{k+2} , \za^0_{k+3}\ldots }]
\eeqn
We want to rewrite this defect word as
$(\za_{z_0-L_0},\ldots,\za_{z_0},\ldots,\za_{z_0+R_0})$, where
$z_0$ is roughly in the center of the defect.  
So let  $L_0:=\lceil w/2 \rceil-1$ and $R_0:=\lfloor w/2 \rfloor$, so
that $w=L_0+R_0+1$.  If $z_0:=i+L_0+1$, then $z_0-L_0=i+1$ and $z_0+R_0=k$, as desired.
Define  $\zd^0_{z}:=\za^0_{z_0+z}$ for all $z\in\CC{-L_0...R_0}$,
and rewrite eqn.(\ref{variable.width.defect.0}) as:
\beqn
\label{variable.width.defect.1}
\ba^0 \ = \
 [\underline{\ldots , \za^0_{z_0 - L_0-2} , \za^0_{z_0 - L_0-1}} \ 
{ \zd^0_{-L_0}...\zd^0_0...\zd^0_{R_0}} \ 
 \underline{\za^0_{z_0 + R_0+1} , \za^0_{z_0 + R_0+2} \ldots }]
\eeqn
[If $w=0$, then $L_0=-1$, \ $R_0=0$, and $z_0=i+1$, so eqn.(\ref{variable.width.defect.1})
is equivalent to the zero-width defect in eqn.(\ref{zero.width.defect}).]
For all $t\in\Natur$, let $\ba^t:=\Phi^t(\ba)$.
We say the defect is {\dfn $\Phi$-persistent} if $\ba^t$ has
a defect for all $t\in\Natur$.  In this case,
\beqn
\label{variable.width.defect}
\ba^t \ = \ 
 [\underline{\ldots ,\za^t_{z_t - L_t-2} , \za^t_{z_t - L_t-1}} \ 
{ \zd^t_{-L_t} ... \zd^t_0 ... \zd^t_{R_t}} \ 
 \underline{ \za^t_{z_t + R_t+1} , \za^t_{z_t + R_t+2} \ldots }]
\eeqn
for some $z_t\in\Zahl$,  $R_t\in\Natur$ and $L_t \in \{R_t, R_t-1\}$.
The next lemma bounds the growth-rate and displacement of the defect
during one $\Phi$-iteration.

\Lemma{\label{instant.velocity}}
{Let $t\in\Natur$.  Then:
\bthmlist
\item $z_t-L_t-1 \leq z_{t+1}-L_{t+1}$ and  $z_{t+1}+R_{t+1} \leq z_t+R_t + 1$.

\item $z_t - L_t - 2 \  \leq \ z_{t+1} \ \leq \ z_t + R_t + 1$.
\ethmlist
}
\bthmprf {\bf(a)} For simplicity, set  $t:=1$.
The boundaries of  
the defect word can advance by at most one unit during each timestep,
because $\Phi$ is a nearest neighbour CA and $\Phi(\bS)\subseteq\bS$.  In other words,
$z_0-L_0-1 \leq z_1-L_1$
and also $z_1+R_{1} \leq z_0+R_0 + 1$.  {\bf(b)} follows because
$z_1-L_1-1\leq z_1 \leq z_1+R_{1}$, because $L_1\geq -1$ and $0\leq R_1$.
\ethmprf

The  width $w_t\approx 2L_t$ of the defect word may
fluctuate with time.  We say that the defect is a {\dfn particle} if
$L:=\D\max_{t\in\Natur}\{L_t\}$
and $R:=\D\max_{t\in\Natur}\{R_t\}$ are finite (possibly $L=-1$ and $R=0$).
Otherwise the defect is called a {\dfn blight} (i.e. its size grows without
bound over time).  We will restrict our attention to  particles.
It will
be convenient to treat the particle as having constant width.  Hence, we
rewrite eqn.(\ref{variable.width.defect}) as
\beqn
\label{variable.width.defect.2}
\begin{array}{rclcccccr}
\ba^t &= &
 [\ldots& \za^t_{z_t - L-2} & \za^t_{z_t - L-1}
  &
{\zd^t_{-L} ... \zd^t_0 ... \zd^t_{R}} &
 \za^t_{z_t + R+1} & \za^t_{z_t + R+2} &\ldots ], \\
 &= &
 [\ldots & \zl^t_{2} & \zl^t_{1}  &
{\zd^t_{-L} ... \zd^t_0 ... \zd^t_{R}} & 
 \zr^t_{1} & \zr^t_{2} &\ldots ], \\
\cline{4-5}\cline{7-8} 
\end{array}
\eeqn
That is: we pad the left side (resp. right side) of the defect with
$L-L_t$ (resp. $R-R_t$) of the `admissible'
symbols, if necessary, and then we define
$\zl^t_n := \za^t_{z_t-L-n}$ and
$\zr^t_n := \za^t_{z_t+R+n}$ for all $n\in\Natur$
(note that, for convenience, we reverse the sign of index $n$ in $\zl^t_n$).
We say $W:=R+L+1$ is the {\dfn width}
of the defect particle.
(If all the defects had zero width, then
$L=-1$ and $R=0$ and $W=0$, so the non-underlined block is empty.)
We can now represent the defect particle as a finite automaton.

A {\dfn finite automaton} is a quintuple
$(\sI,\sD,\sO;\Ups,\Omg)$, where $\sI$ is a finite {\dfn input
alphabet}, $\sD$ is a finite {\dfn state domain}, $\sO$ is a finite
{\dfn output alphabet}, $\Ups:\sI\x\sD\into\sD$ is an {\dfn update
rule}, and $\Omg:\sI\x\sD\into\sO$ is an {\dfn output rule}.  Finite
automata are models of simple computers: starting from initial state
description $\zd_0\in\sD$,
 and fed an input stream $(\zi_0,\zi_1,\zi_2,\ldots)\in\sI^\Natur$,
the automaton undergoes a series of state transitions $\zd_0\leadsto
\zd_1\leadsto \zd_2 \leadsto \ldots$ [where $\zd_{t+1}:=\Ups(\zi_t,\zd_t)$] and
produces an output stream $(\zo_1,\zo_2,\zo_3,\ldots)\in\sO^\Natur$, where
$\zo_{t+1}:=\Omg(\zi_t,\zd_t)$.  See  \cite[\S2.2]{HopcroftUllman}.

 The defect particle in eqn.(\ref{variable.width.defect.2}) behaves like a finite automaton  with $\sI:=\sA^{L+2}\x\sA^{R+1}$, \  $\sD:=\sA^\CC{-L...R}$,
 and $\sO:=\CC{-L\!-\!2\ldots R\!+\!1}$. 
  The automaton's inputs are $\bl^t:=[\zl_{L+2},\ldots, \zl^t_1]$
and $\br^t:=[\zr_1,\ldots, \zr_{R+1}]$, its internal state is 
 $\bd^t:=[\zd^t_{-L}, \ldots, \zd^t_0, \ldots, \zd^t_{R}]\in\sD$,
and its output is a `velocity vector' in $\dV:=\CC{-L\!-\!2\ldots R\!+\!1}$.
That is, there is a unique {update rule}
$\Upsilon:\sA^{L+2}\x\sD\x\sA^{R+1} \into \sD$ and {\dfn velocity function}
$\vV:\sA^{L+2}\x\sD\x\sA^{R+1} \into \dV$ such that
\beqn
\label{defect.finite.state.update}
\bd^{t+1}\quad=\quad  \Upsilon(\bl^t,\bd^t,\br^t)
\quad\And\quad
z_{t+1}= z_t + \vV(\bl^t,\bd^t,  \br^t) 
\eeqn
Let $\gL$ and $\gR$ be the unique $(\shift{},\Phi)$-transitive
components of $\bS$ such that $[\ldots \zl^t_3 \zl^t_2 \zl^t_1]$ is
$\gL$-admissible and $[\zr^t_1 \zr^t_2 \zr^t_3 \ldots ]$ is $\gR$-admissible
for all $t>0$ (possibly $\gL=\gR$).  We say that $\ba^t$ has an {\dfn
$(\gL,\gR)$ defect particle of width $W$}.

\begin{figure}[h]
\centerline{\footnotesize
\begin{tabular}{|c|c|c|c|}
\hline
\includegraphics[scale=6,trim=0 2 0 0 ]{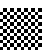} &
\includegraphics[scale=6,trim=0 2 0 0 ]{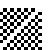} &
\includegraphics[scale=6,trim=0 2 0 0 ]{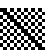} &
\includegraphics[scale=6,trim=0 2 0 0 ]{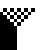} \\
$(*)$ & $\boldsymbol{(\gamma^-)}$ & $\boldsymbol{(\gamma^+)}$ &
$\boldsymbol{(\beta)}$   \\
\hline\hline
\includegraphics[scale=6,trim=0 2 0 0 ]{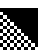} &
\includegraphics[scale=6,trim=0 2 0 0 ]{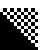} &
\includegraphics[scale=6,trim=0 2 0 0 ]{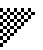} &
\includegraphics[scale=6,trim=0 2 0 0 ]{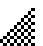} \\
$\boldsymbol{(\alpha^+)}$ & $\boldsymbol{(\omg^+)}$ &
 $\boldsymbol{(\alp^-)}$ & $\boldsymbol{(\omg^-)}$ \\
\hline
\end{tabular} }
\caption{{\footnotesize 
$(*)$  The periodic background generated
by $\ECA{184}$ acting on  $\bG_*$
 in Example \ref{X:variable.width.defect}(a).\quad
$\boldsymbol{(\alpha^\pm,\beta,\omg^\pm,\gamma^\pm)}$:
Defect particles of  $\ECA{184}$ acting on $\bG$.
(Note that the image labelled $\boldsymbol{(\bet)}$ actually depicts
the coalescence of an $\omg^+$ and an $\alp^-$ into  a $\bet$.)
See also \cite[Examples 2.2(a) \& 3.5(c)]{PivatoDefect1}
or  \cite[\S III(A)]{BNR91}.
\label{fig:184.defects}}}
\end{figure}

\example{\label{X:variable.width.defect}
(a) (ECA\#184) Let $\sA=\{0,1\}$. 
Let $\ECA{184}:\AZ\into\AZ$ be elementary cellular
automaton \#184. (Recall: the number `184' encodes the local rule
$\phi:\sA^{\{-1,0,1\}}\into\sA$ via the formula
$\sum_{i=0}^1 \sum_{j=0}^1 \sum_{k=0}^1 \phi(i,j,k) (4i + 2j + k) \ = \ 184$).
Let $\bG_*\subset\AZ$ be
the Markov subshift given by digraph $\black\leftrightarrows \white$
(we use the convention that $\black=0$ and $\white=1$, to ease comparison
between equations and figures).
Thus $\bG_*:=\{(\black\white)^\oo, \ (\white\black)^\oo\}$,
where $(\black\white)^\oo:=[\ldots\black\white\underline{\black}\white\black\white\ldots]$, etc.
(the zeroth coordinate is underlined). 
Then $\ECA{184}\restr{\bG_*} = \shift{}$,
as shown in Figure \ref{fig:184.defects}$(*)$.  There are
two $(\bG_*,\bG_*)$-defects of width 0, shown in 
Figure \ref{fig:184.defects}$(\gam^\pm)$.  The $\gam^+$ defect
 consists of a single inadmissible transition $(\black,\black)$, while  
$\gam^-$ consists of an inadmissible $(\white,\white)$ transition.
Now, $W=0$, so $\sD=\emptyset$ and
$\Ups$ is trivial, and $\vV:\sA\x\sA\into\{-1,0,1\}$.
We have $\vV(0,0)=1$ (for $\gam^+$)  and $\vV(1,1)=-1$ (for $\gam^-$).

Let $\bG_{01}\subset\AZ$ be the Markov subshift determined by digraph
$\stackrel{\curvearrowleft}{\black} \quad \stackrel{\curvearrowleft}{\white}$.  Thus,
$\bG_{01}=\bG_0\disj\bG_1$, where  
$\bG_0:=\{\black^\oo\}$ and $\bG_1:=\{\white^\oo\}$.  Note that
$\bG_{01} \subset \Fix{\ECA{184}}$; hence $\bG_{01}$ is a $\Phi$-invariant
subshift.
The $(\bG_0,\bG_1)$-defect
of width $0$ is a particle shown in Figure \ref{fig:184.defects}$(\bet)$.
(The $(\bG_1,\bG_0)$-defect of width zero is unstable, and immediately
`decays' into a two other particles).
Again, $W=0$, so $\sD=\emptyset$,
$\Ups$ is trivial, and $\vV:\sA\x\sA\into\{-1,0,1\}$.
We have $\vV(0,1)=0$;  i.e. the $\beta$ particle is stationary.

Let
$\bG\subset\AZ$ be the subshift of finite type determined by the admissible
$3$-tuples $\bG_3:=\{(\black\black\black),(\white\white\white),(\white\black\white),(\black\white\black)\}$.  We block-recode this
as a Markov subshift in the alphabet $\sA^3$, given by digraph
$\stackrel{\curvearrowleft}{(\black\black\black)}; \quad (\white\black\white)\leftrightarrows(\black\white\black);
\quad \stackrel{\curvearrowleft}{(\white\white\white)}$.  Thus,
$\bG=\bG_{*}\disj\bG_{01}$, where  
$\bG_{*}$ and $\bG_{01}$ are as above.
Thus $\bG$ has three $(\ECA{184},\shift{})$-transitive components:
$\bG_*$, $\bG_0$ and $\bG_1$.
Figure \ref{fig:184.defects} shows several defect particles of
$\ECA{184}$ acting on $\bG$.
 The defects $\alp^\pm$ and $\omg^\pm$ have width 1, so
$\sD=\sA^3$.   Although they had width 0 as
defect particles in $\bG_*$ or $\bG_{01}$, 
the defects $\gam^\pm$ and $\bet$ have width 2 as defect particles in $\bG$,
so $\sD=\sA^3\x\sA^3$ (although,
by the definition of block-recoding, we could replace this with $\sD=\sA^4$).
For all seven particles, the defect word $\bd\in\sD$ is constant
over time.
The values of $\gL$, $\gR$, $\bd$, and the (constant) value of 
$\vV:\sA\x\sD\x\sA\into\{-1,0,1\}$
for each defect are given by Table \ref{table:184.defects}.
See also Example \ref{X:ballistic}(a).

\begin{table}
\[\begin{array}{r|c|c|c|c|c|c|c|c|}
\cline{2-8}
&\multicolumn{4}{|c|}{W=1} &
\multicolumn{3}{|c|}{W=2} \\
\cline{2-8}
& \alp^- & \alp^+ &\omg^- & \omg^+ & \gam^{-} & \gam^+ &\bet  \\
\hline  
\gL            & \bG_* & \bG_*  & \bG_1 & \bG_0  &  \bG_*  &\bG_*  &  \bG_0  \\ 
\gR            & \bG_1 & \bG_0  & \bG_* & \bG_*  &  \bG_*  & \bG_*  &  \bG_1  \\ 
\bd &  \black\white\white&  \white\black\black & \white\white\black& \black\black\white& \black\white\white\black& \white\black\black\white & \black\black\white\white\\
\vV & -1& 1 & -1 & 1 & -1 & 1 & 0\\ 
\hline
\end{array}
\]
\caption{\footnotesize Seven defect particles for ECA\#184 acting on $\bG$; see
Example \ref{X:variable.width.defect}(a).\label{table:184.defects}}
\end{table}

{
\begin{figure}
\centerline{
\begin{tabular}{|c|c|c|c|c|}
\hline
\includegraphics[scale=4.7]{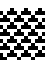} &
\includegraphics[scale=4.7]{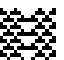} &
\includegraphics[scale=4.7]{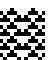} &
\includegraphics[scale=4.7]{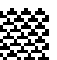} &
\includegraphics[scale=4.7]{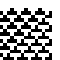} \\
\hline
$(*)$  & $\boldsymbol{(\alpha)}$ & $\boldsymbol{(\beta)}$ &  $\boldsymbol{(\gamma^+)}$ & $\boldsymbol{(\gamma^-)}$  \\
\hline 
\end{tabular} }
\caption{{\footnotesize 
$(*)$  The periodic background generated by $\ECA{54}$
acting on $\gB$ in Example \ref{X:variable.width.defect}(b).\quad
$\boldsymbol{(\alpha, \beta,\gamma^\pm)}$: Four defect particles
of $\ECA{54}$ acting on $\gB$.
See also \cite[Example 3.5(b)]{PivatoDefect1},
\cite[Fig.8]{CrHa97}, or \cite[\S III(C)]{BNR91}.
\label{fig:54.defects}}}
\end{figure}
}

(b) (ECA\#54)  Let $\sA=\{0,1\}$.
Let $\ECA{154}:\AZ\into\AZ$ be elementary cellular
automaton \#54.
Let $\gB:=\gB_0\disj\gB_1$,
where $\gB_0$ is the four-element $\shift{}$-orbit of $\overline{0010}$ and 
$\gB_1$ is the  four-element $\shift{}$-orbit of $\overline{1101}$.
Then $\ECA{54}(\gB_0)=\gB_1$ and
$\ECA{54}(\gB_1)=\gB_0$, so $\gB$ is $(\ECA{54},\shift{})$-transitive,
 so all defects have $\gL=\gB=\gR$.
Also, $\ECA{54}^2\restr{\gB}=\shift{2}$ 
[see Figure \ref{fig:54.defects}$(*)$].
We recode $\gB$ as a topological Markov subshift  in the alphabet
$\sA^4$, with admissible $4$-words
$\sB:=\{\white\white\white\black,$ $ \  
\white\white\black\white,$ $ \  
\white\black\white\white,$ $ \  
\black\white\white\white; \ 
\black\black\black\white,$ $ \ 
\black\black\white\black,$ $ \ 
\black\white\black\black,$ $ \ 
\white\black\black\black\}$.
Figure \ref{fig:54.dislocate} shows the $\ECA{54}$-evolution of
several defect particles in $\gB$,
along with the relevant values of  $z$, $R$, $L$, $\vV$,
$\sD$, and $\Ups$. See also Example \ref{X:ballistic}(b).

\begin{figure}[h]
\centerline{\footnotesize
\psfrag{D5}[][]{$\sD=\sA^\CC{-2...2}\cong \sA^5$}
\psfrag{D4}[][]{$\sD=\sA^\CC{-1...2}\cong \sA^4$}
\psfrag{D1}[][]{$\sD=\sA$}
\psfrag{Y}[][]{$\Upsilon$}
\psfrag{l}[][]{$L$}
\psfrag{r}[][]{$R$}
\psfrag{z}[][]{$z$}
\psfrag{v}[][]{$\vV$}
\psfrag{a0}[][]{$\alp_0$}
\psfrag{a1}[][]{$\alp_1$}
\psfrag{a2}[][]{$\alp_2$}
\psfrag{a3}[][]{$\alp_3$}
\psfrag{b0}[][]{$\bet_0$}
\psfrag{b1}[][]{$\bet_1$}
\psfrag{b2}[][]{$\bet_2$}
\psfrag{b3}[][]{$\bet_3$}
\psfrag{g0}[][]{$\gam^+_0$}
\psfrag{g1}[][]{$\gam^+_1$}
\psfrag{G0}[][]{$\gam^-_0$}
\psfrag{G1}[][]{$\gam^-_1$}
\includegraphics[scale=0.63,angle=-90]{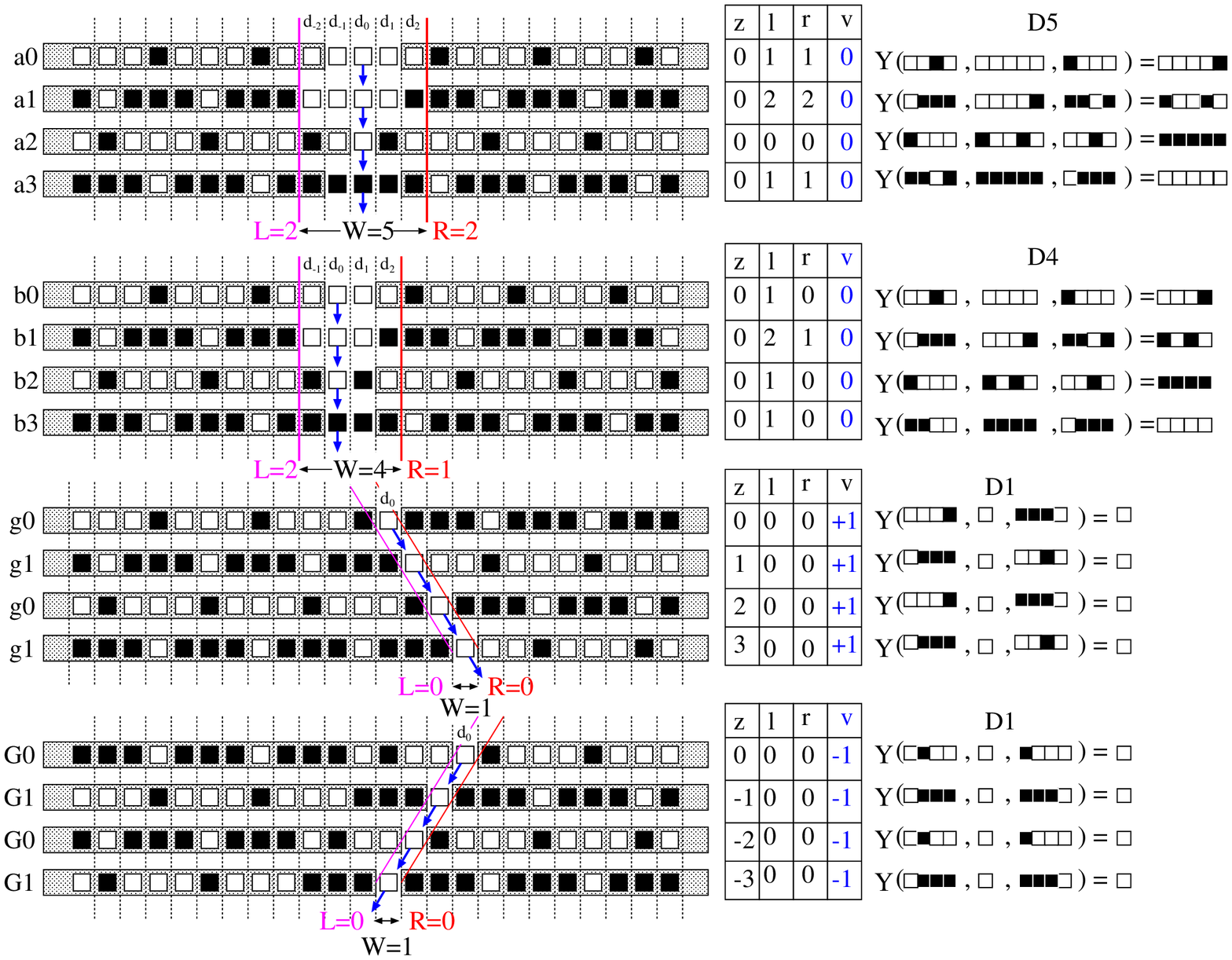}}
\caption{{\footnotesize  Defects in ECA\#54.
We treat the symbol $\za_z$ as `defective' if the word
$(\za_{z-2},\za_{z-1},\za_{z},\za_{z+1})$ is not $\gB$-admissible.
The admissible segments of each sequence are boxed; hence the unboxed
segments are the defect words.  The table on the right describes the
values for $z$, $R$, $L$, and $\vV$ and
the definition of $\sD$ in each case, as well as the relevant
values of the update rule $\Upsilon:\sA\x\sD\x\sA\into\sD$.
See Example \ref{X:variable.width.defect}(b).
\label{fig:54.dislocate}}}
\end{figure}

\begin{figure}
\centerline{
\begin{tabular}{|c|c|c|}
\hline
\includegraphics[height=9em,width=10.2em,trim=0 2 0 0 ]{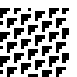} &
\includegraphics[height=9em,width=10.2em,trim=0 2 0 0 ]{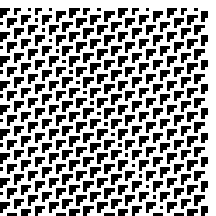} &
\includegraphics[height=9em,width=10.2em,trim=0 2 0 0 ]{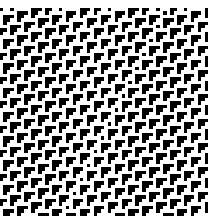} \\
$(*)$     & {\bf(A)}   & {\bf(B)} \\
\hline
\end{tabular} }
\caption{{\footnotesize $(*)$  A $30\x30$ image of
the periodic spacetime diagram of $\ECA{110}$ acting on $\gE$
from Example \ref{X:variable.width.defect}(c);\quad
{\bf(A,B)}  $60\x60$ images of the
$\ECA{110}$-evolution of two defect particles in $\gE$;
See \cite{Lind2,McIntosh1,McIntosh2,CrHR}, \cite[\S3.1.4.4]{Ilachinski},
 \cite[Chap.11]{Wolfram2}, \cite[Example 3.5(d)]{PivatoDefect1},
 and especially \cite{Cook}.
\label{fig:110.defects}}}
\end{figure}

\begin{figure}
\centerline{
\psfrag{DA}[][]{$\sD=\sA^\CC{-5...5}\cong \sA^{11}$}
\psfrag{DB}[][]{$\sD=\sA^\CC{-6...6}\cong \sA^{13}$}
\psfrag{Y}[][]{$\Upsilon$}
\psfrag{l}[][]{$L$}
\psfrag{r}[][]{$R$}
\psfrag{z}[][]{$z$}
\psfrag{v}[][]{$\vV$}
\includegraphics[angle=-90,scale=0.6]{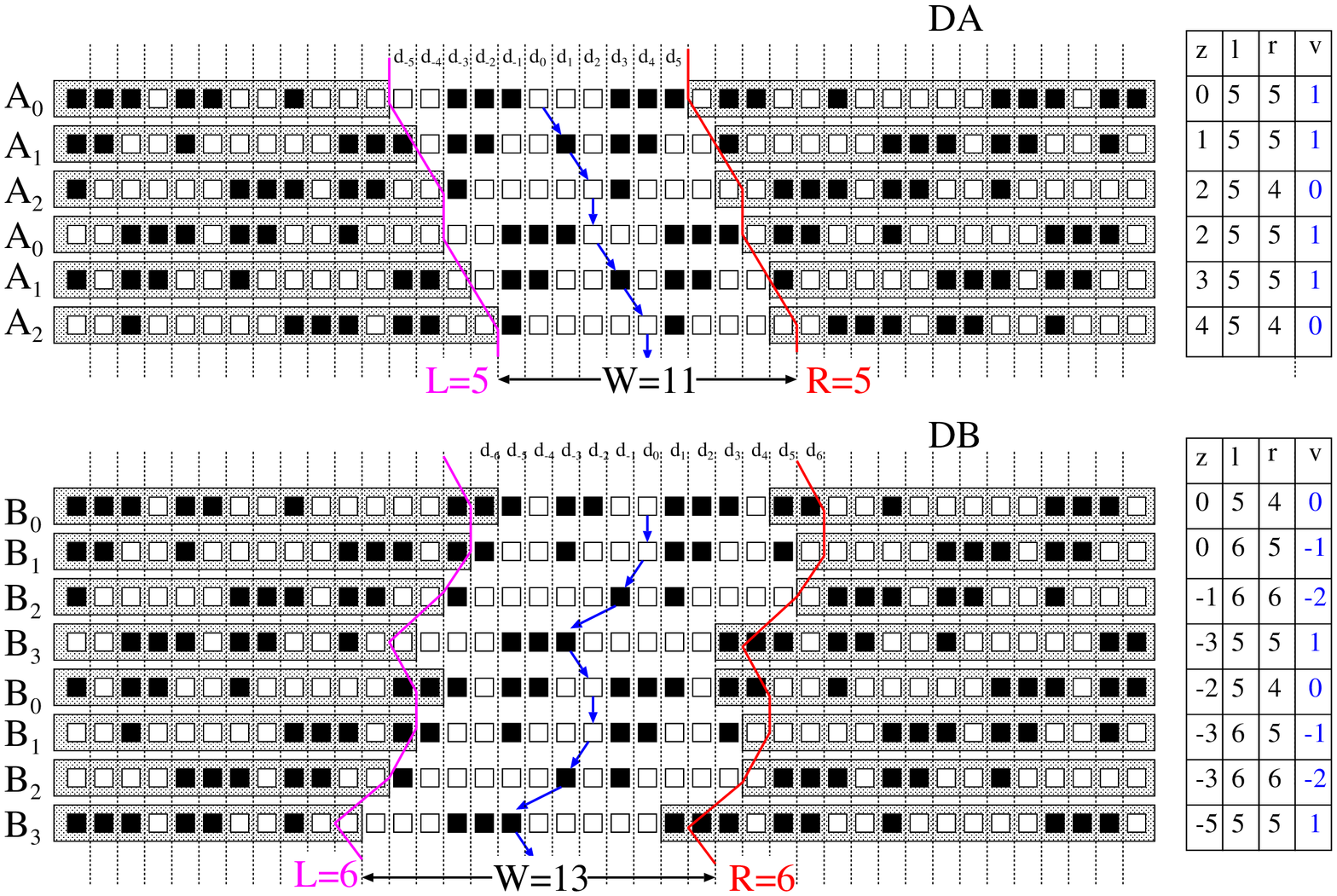}}
\caption{{\footnotesize  The $A$ and $B$ defect particles of
ECA\#110.
We treat the symbol $\za_z$ as `defective' if the word
$(\za_{z-6},\ldots,\za_{z},\za_{z+1},\ldots,\za_{z+7})$ is not $\gE$-admissible.
The admissible segments of each sequence are boxed; hence the unboxed
segments are the defect words.  The table on the right describes the
values for $z$, $R$, $L$, and $\vV$ in each case.
The arrow path is the sequence $(z_t)_{t\in\Natur}$.
The left-hand and right-hand polygonal paths are the sequences 
$(L_t)_{t\in\Natur}$ and  $(R_t)_{t\in\Natur}$.
See Example \ref{X:variable.width.defect}(c).
\label{fig:110.prop}}}
\end{figure}

(c) (ECA\#110) Let $\sA=\{0,1\}$.
Let $\gE\subset\AZ$ be
the 14-element $\shift{}$-orbit of the 14-periodic sequence
%$\overline{00010011011111}$.  
$(\black\black\black\white\black\black\white\white\black\white\white\white\white\white)^\oo$.  
If $\ECA{110}$ is ECA\#110, then $\ECA{110}\restr{\gE}=\shift{4}$
[see Figure \ref{fig:110.defects}$(*)$], so
$\gE$ is $(\ECA{110},\shift{})$-transitive, so all defects have
$\gL=\gE=\gR$.
 Figure \ref{fig:110.prop} shows the $\ECA{110}$-evolution
of two defects (called  `$A$' and
`$B$' in the nomenclature of \cite{Cook})
along with the relevant values of  $z$, $R$, $L$, $\vV$,
and $\sD$.
}

\NRemark{\label{defect.particle.remarks}
(a) \ $\vV$ takes values in $\dV:=\CC{-L\!-\!2\ldots R\!+\!1}$ by
Lemma \ref{instant.velocity}(b).  However, the average value of
$\vV$ over time must be in $\CC{-1,1}$, because 
the left endpoint of the defect has a minimum velocity of $-1$,
while the right endpoint has a maximum velocity of $+1$
[by Lemma \ref{instant.velocity}(a)].
If $\vV<1$ (resp. $\vV>1$), this means that the right (resp. left)
endpoint is moving  leftward (resp. rightward) at speed greater
than 1, which means the defect particle is shrinking, which is only sustainable
for a short period of time.  For example, the particle can achieve
an instantaneous velocity $\vV=R+1$ only by shrinking from a defect
of width $W$ to one of width $0$;  it must later remain at velocity $\vV=0$ 
for $(R+1)$ iterations to grow back to width $W$.  

 The `constant width' convention of eqn.(\ref{variable.width.defect.2})
masks this  shrinkage by `padding' the defect word $\bd^{t+1}$
with up to $R+1$ admissible characters from 
$\ba^{t+1}_{\CC{z_{t+1}\ldots z_{t+1}+R+2}}$;  this is why
the function $\Ups$ needs $\br^t$ as 
input.  Likewise, possibly rapid leftward motion requires $\Ups$ to incorporate
$\bl^t$ as input.  In most examples,
however, the particle moves slowly, 
and we can reduce the number of boundary inputs.

(b)  If $W=0$, then  $\sD=\emptyset$ and $\Ups$ is trivial, while
$\vV$ is a function $\vV:\sA\x\sA\into\{-1,0,1\}$.

(c)  If $W\geq 1$, then 
  by passing to the $W$th higher power representation \cite[Defn.1.4.4]{LindMarcus}, we can assume that $L=0$ and $R=1$, so that $W=2$. 
To see this, replace
$\sA$ with $\hsA:=\sA^W$, and represent 
$\ba=[\ldots \za_{-1} \ \za_0 \ \za_1 \ \za_2 \ldots]\in\AZ$ by
\[
\hba\ := \
 \lb[\mbox{\footnotesize $\Array{\ldots\\ \\ \ldots} \  \Array{ \za_{-2W} \\ \vdots \\ \za_{-W-1}}  \ 
\Array{ \za_{-W} \\\vdots \\ \za_{-1}}  \
\Array{\za_{0}\\ \vdots \\ \za_{W-1}} \
\Array{\za_{W} \\ \vdots \\ \za_{2W-1} }
\ \Array{ \za_{2W} \\ \vdots \\ \za_{3W-1}} \ \Array{\ldots\\ \\ \ldots} $ }\rb]
\ = \
[\ldots \hza_{-1} \ \hza_0 \ \hza_1 \ldots] \ \in \ \hsA^\Zahl
\]
(Note:  this is {\em not} the same as the higher {\em block} recoding
described earlier). 
Thus, if $\ba^t$ is as in eqn.(\ref{variable.width.defect.2}), 
and  $\hz_t := \lfloor z_t/W \rfloor$, then
%\beqn\label{variable.width.defect.3}\eeqn
$\hba^t \ = \ [\ldots \hzl_{2} \ \hzl_{1} \ \h\zd_{0} \ \h\zd_1 \
\hzr_{1} \ \hzr_{2} \ldots]$, where $\hzl_n,\hzr_n \in \bS_{W}$ for all
$n\in\Natur$, while $\hzd_0$ and $\hzd_1$ are in $\hsA$.  The original defect
word $\bd^t$ is split between $\hzd_{0}$ and $\hzd_1$.  The particle's
behaviour now depends only on its nearest neighbours, and its speed is
never greater than 1.  In other words, we have
$\hUps:\hsA\x\sD\x\hsA\into\sD$ and
$\vV:\hsA\x\sD\x\hsA\into\{-1,0,1\}$ in
eqn.(\ref{defect.finite.state.update}).  The price
of this manoeuvre is that each defect
word $\bd\in\sD$ can be represented in $W$ distinct ways as a pair
$(\h\zd_{-1},\h\zd_0)$, depending upon the value of $z_t$ mod $W$; this
may translate into $W$ spuriously distinct `particle types' [see
Definition \ref{def:type} and Remark \ref{rem:defect.random.walk}(a)
below].  Also, it may make some dislocations look like interfaces [see
Remark (d) below].  Nevertheless, it will be useful in the
proofs of Proposition \ref{Turing.CA}, Theorem
\ref{defect.random.walk} and Corollary \ref{defect.random.walk2}.

(d) Two sequences $\bb,\bc\in\AZ$ are {\dfn homoclinic} if there is
some $N>0$ such that $\zb_z = \zc_z$ for all $z\in\Zahl$ with $|z|>N$.
The defect in $\ba^t$ is called {\dfn removable} if $\ba^t$ is
homoclinic to some $\bs\in\bS$; see \cite[\S1]{PivatoDefect1}.
Otherwise the defect in $\ba^t$ is {\dfn essential} ---i.e. it is
impossible to remove the defect by changing $\ba^t$ in some finite
set.  If $\gL\neq\gR$ [e.g. Figure \ref{fig:184.defects}$(\beta)$]
then an $(\gL,\gR)$-defect is called an {\dfn interface}, and is
necessarily an essential defect; see \cite[\S2]{PivatoDefect1}.  If
$\gL=\gR$ [e.g. Figure \ref{fig:184.defects}$(\gam^\pm)$], then an
$(\gL,\gR)$-defect is called a {\dfn dislocation}, and may or may not
be essential, depending upon whether it induces a `phase slip' in the
periodic structure of $\gL$; see \cite[Example 3.1]{PivatoDefect1}.  }

\begin{table}
{\footnotesize
\centerline{
\psfrag{S1}[][]{{\scriptsize (\S\ref{S:ballistic})}}
\psfrag{S3}[][]{{\scriptsize (\S\ref{S:turing})}}
\psfrag{S2}[][]{{\scriptsize (\S\ref{S:diffuse})}}
\psfrag{C2.1}[][]{{\tiny (Theorem \ref{defect.random.walk})}}
\psfrag{C2.2}[][]{{\tiny (Corollary \ref{defect.random.walk2})}}
\includegraphics[scale=0.7]{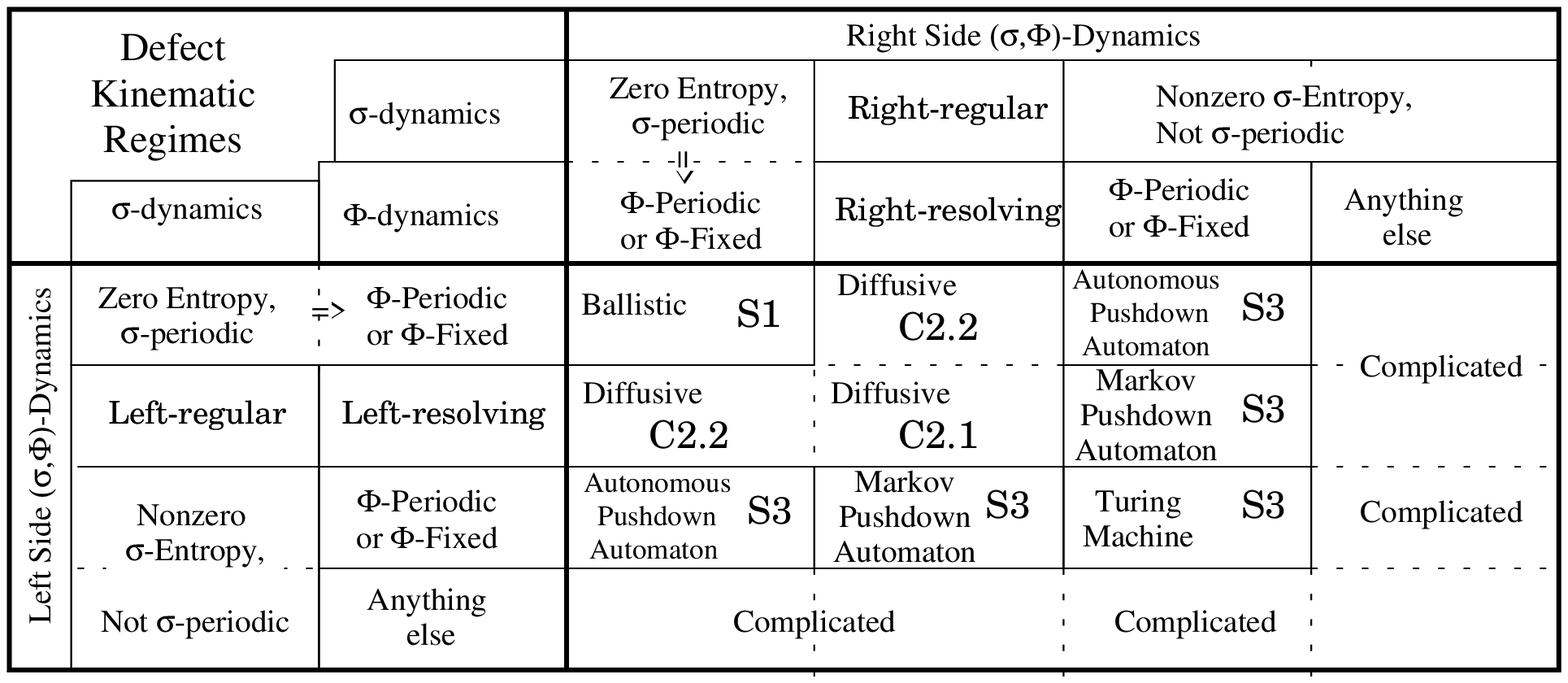}}
\caption{Kinematic regimes for one-dimensional defect particles.
\label{table:defect.regime}}}
\end{table}

\begin{sloppypar}
The kinematics of a one-dimensional defect particle falls into several
regimes summarized in Table \ref{table:defect.regime}, depending on
the $(\shift{},\Phi)$-dynamical complexity of $\gL$ and $\gR$.  In the
{\em Ballistic} regime (see \S\ref{S:ballistic}), the defect acts as a
finite automaton driven by periodic input, and moves with a constant
average velocity through a periodic background.  ECAs \#54, \#62,
\#110, and \#184 are all in this regime, which has been studied
empirically in \cite{Gra83,Gra84,BoRo91,BNR91,Han,CrHa97,CrHM}.  At
the opposite extreme, in the {\em Diffusive} regime (see
\S\ref{S:diffuse}), the defect acts like a finite-state Markov
process, and performs a generalized random walk.  Diffusive defect
dynamics has previously been analyzed by Eloranta
\cite{ElNu,Elo93a,Elo93b,Elo94}.  In the {\em Turing} regime (see
\S\ref{S:turing}), the defect moves through an inert, positive-entropy
background, and modifies this background with its passing; the system
acts like a Turing machine, where the particle is the `head' and the
inert background is the `tape'.  In the {\em Autonomous Pushdown
Automaton} regime (see \S\ref{S:turing}), the defect has a
$\Phi$-fixed, positive $\shift{}$-entropy domain on one side (which we
treat as a `stack' memory), and a zero-entropy domain on the other
side; the system acts like a pushdown automaton operating autonomously
(i.e. without external input).  In the {\em Markov Pushdown Automaton}
regime (see \S\ref{S:turing}), the defect has a $\Phi$-fixed, positive
$\shift{}$-entropy domain on one side (which we treat as a `stack'),
and a $\Phi$-resolving subshift on the other; the system acts like a
pushdown automaton driven by noise from a Markov process.  The {\em
Complicated} regime is none of the above, and is probably too diverse
to make any useful generalizations.
\end{sloppypar}

\section{The Ballistic Regime \label{S:ballistic}}

Let $\Phi:\AZ\into\AZ$ be a cellular automaton, and 
let $\bX\subset\AZ$ be a $\Phi$-invariant Markov shift.
Let $\bL,\bR\subset\bX$ be $(\Phi,\shift{}$)-transitive subshifts of $\bX$,
let $W\in\Natur$,
and let $\bD^W_{\bL,\bR}$ be the set of all sequences in $\AZ$
with a single $(\bL,\bR)$ defect particle of width $W$,
such as shown in eqn.(\ref{variable.width.defect.2}). 
By hypothesis, $\Phi(\bD^W_{\bL,\bR})\subseteq \bD^W_{\bL,\bR}$.

\Theorem{\label{thm:ballistic}}
{
  Suppose $\bL$ and $\bR$ are $\shift{}$-periodic and  $(\shift{},\Phi)$-transitive.
Then the dynamical system $(\bD^W_{\bL,\bR},\Phi)$ is isomorphic
to a dynamical system $(\sX\x\Zahl,\Xi)$, where $\sX$
is a finite set, and where $\Xi:\sX\x\Zahl\into\sX\x\Zahl$ is
defined by $\Xi(\zx,z):=(\xi(\zx), z+\vV(\zx))$ for some
functions $\xi:\sX\into\sX$ and $\vV:\sX\into\dV:=\CC{-L\!-\!2\ldots R\!+\!1}$.
}
{\em Proof idea:} 
The defect particle's internal state is a finite automaton driven by a periodic
input (because $\bL$ and $\bR$ are periodic);
thus, by incorporating the phase of this periodic input into the state
description of the defect, we can treat it as an autonomous
finite automaton (i.e. a finite-state dynamical system) $(\sX,\xi)$.
The defect's position 
is then obtained by integrating the velocity signal generated
by  $(\sX,\xi)$.

\bthmprf  Any $\sigma$-periodic sequence is automatically $\Phi$-periodic.
Thus, by hypothesis, every element of $\bL$ is $\shift{P_L}$-fixed
and $\Phi^{Q_L}$-fixed for some $P_L,Q_L\in\Natur$.  But $\bL$
is $(\Phi,\shift{})$-transitive, so this means that
$\bL$ consists of a single finite $(\shift{},\Phi)$-orbit
containing exactly $P_L Q_L$ elements. 
Recall that $\bL^-\subset\sA^\OC{-\oo...-1}$ is the set of all left-infinite
$\bL$-admissible sequences.  Since
$\bL$ is a $\shift{}$-periodic Markov shift, any element
$[\ldots,\zl_{3},\zl_{2},\zl_{1}]$ in $\bL^-$ is
completely determined by the entry $\zl_{1}$.  Hence
there is a subset $\sL\subseteq\sA$ with $\#(\sL)=
P_L Q_L$, and a bijection $\lam:\sL\into\bL^-$
where, for any $\zl\in\sL$, $\lam(\zl)$ is the unique
sequence $[\ldots,\zl_{3},\zl_{2},\zl_{1}]$ in $\bL^-$ 
with $\zl_{1}=\zl$.  Furthermore, there are
bijections $\varphi_L:\sL\into\sL$ and
$\varsigma_L:\sL\into\sL$ such that $\Phi\circ\lam=\lam\circ\varphi_L$
and $\sigma\circ\lam=\lam\circ\varsigma_L$.

Likewise,  $\bR$ is $\shift{P_R}$-fixed
and $\Phi^{Q_R}$-fixed (for some $P_R,Q_R\in\Natur$), so
$\bR$ has exactly $P_R Q_R$ elements. 
Recall that $\bR^+\subset\sA^\CO{1...\oo}$ is the set of right-infinite
$\bR$-admissible sequences.
There is a subset $\sR\subseteq\sA$ with
$\#(\sR)=P_R Q_R$, and a bijection $\rho:\sR\into\bR^+$
so that, for any $\zr\in\sR$, $\rho(\zr)$ is the unique
sequence $[\zr_1,\zr_2,\zr_3,\ldots]$ in $\bR^+$ with with $\zr_{1}=\zr$.
There are  bijections $\varphi_R:\sR\into\sR$ and
$\varsigma_R:\sR\into\sR$ such that $\Phi\circ\rho=\rho\circ\varphi_R$
and $\sigma\circ\rho=\rho\circ\varsigma_R$.

\begin{sloppypar}
Thus, the sequence $\ba^t$ in eqn.(\ref{variable.width.defect.2}) is
entirely determined by the data $(\zl^t_{1}, \bd^t,
\zr^t_{1}; z_t)\in\sL\x\sD\x\sR\x\Zahl$, because
$[\ldots \zl^t_3, \zl^t_2, \zl^t_1]
  =   \lam(\zl^t_1)$, and
$\rho(\zr^t_1) =  [\zr^t_1, \zr^t_2, \zr^t_3,\ldots ]$.
Define $\Psi:\sL\x\sD\x\sR\x\Zahl\into\bD^W_{\bL,\bR}$
by
$\Psi(\zl,\bd,\zr;z) := [\bl \ \bd \  \br]$, where $\bl:=\lam(\zl)\in\bL^-$,
$\br:=\rho(\zr)\in\bR^+$, and we place $\bd$ so that its center
coordinate is at $z$.  Then $\Psi$ is a bijection.
\end{sloppypar}

Let $\sX:=\sL\x\sD\x\sR$.
If $\Upsilon$ and $\vV$ are as in  eqn.(\ref{defect.finite.state.update}),
then we can restrict them to functions
$\Upsilon_|:\sX\into\sD$ and $\vv:=\vV_|:\sX\into\dV$.
Define $\xi:\sX\into\sX$ by $\xi(\zl, \bd, \zr) \ := \ (\zl', \bd', r')$,
where 
\[
\zl' \ := \ \varsigma_L^v  \circ \varphi_L(\zl),\ \ 
\bd' \ := \  \Ups(\zl,\bd,\zr),
 \ \And \ 
r' \ := \ \varsigma_R^v  \circ \varphi_R(\zr),\
\mbox{where $v:=\vV(\zl,\bd,\zr)$.}
\]
Now define $\Xi:\sX\x\Zahl\into\sX\x\Zahl$ as in the theorem statement.
Then $\Phi\circ\Psi = \Psi\circ\Xi$. 
\ethmprf

\Definition{\label{def:type}}
{The (finite) dynamical system $(\sX,\xi)$ decomposes into a finite
disjoint union of finite $\xi$-orbits, called {\dfn particle types}.
If $\sP\subset\sX$
is a particle type, then $P:=\#(\sP)$ is the {\dfn period} of type $\sP$, and
$\vV(\sP):=\D \frac{1}{P}\sum_{\zp\in\sP} \vV(\zp)$ is the
{\dfn average velocity} of type $\sP$.}

\example{\label{X:ballistic}
(a) (ECA\#184) We continue Example \ref{X:variable.width.defect}(a).
  There are three $(\shift{},\ECA{184})$-transitive components
in $\gG$, so there are three possible choices for $\bR$; for each one,
we list the corresponding values of $\sR$, $P_R$, $Q_R$, $\rho$,
$\varsigma_R$, and $\varphi_R$ in Table \ref{fig:ballistic.184}(A) (the values
for $\sL$, $P_L$, $Q_L$, $\lam$,
$\varsigma_L$, and $\varphi_L$ would be exactly the same).
In Example \ref{X:variable.width.defect}(a), we introduced seven defect
particles for $(\gG,\Phi)$:  four of width 1, and three of width 2.
In all seven cases, we have $\Ups(\zl,\bd,\zr)=\bd$.
Thus,  $\xi:\sL\x\sD\x\sR\into\sL\x\sD\x\sR$ is given by
$\xi(\zl,\bd,\zr)=  \lb(\varsigma_L^v  \circ \varphi_L(\zl), \ \bd,
\varsigma_R^v  \circ \varphi_R(\zr)\rb)$, where $v:=\vV(\zl,\bd,\zr)\in\{-1,0,1\}$.
The value of  $\vV$ is constant
for each particle type, and was shown in the bottom row of 
Table \ref{table:184.defects}.
In all cases, we end up with $\xi=\Id{}$, so all particle types have period 1.
Hence, the average velocity of each type is just the value of
$\vV$ on the (unique) member of that type.

\begin{table}
{\footnotesize
\[
\begin{array}{|r|c|c|c|c|c|c|}
\hline
\bR  & P_R & Q_R & \sR  & \rho:\sR\into\bR & \varsigma_R &\varphi_R \\
\hline			 
\gG_0 & 1 & 1 &   \{0\} &\rho(0)=(000...) & \Id{} & \Id{}\\
\hline
\gG_1 & 1 & 1 &   \{1\} &\rho(1)=(111...) & \Id{} & \Id{} \\ 
\hline
\gG_* & 2 & 2 &   \{0,1\} &{{\rho(0)=(0101...)}\atop{\rho(1)=(1010...)}} &
{{\varsigma_R(0)=1}\atop{\varsigma_R(1)=0}} &
{{\varphi_R(0)=1}\atop{\varphi_R(1)=0}}  \\
\hline
\end{array}
\]
\caption{{\footnotesize Ballistic defects in ECA\#184; 
See Example \ref{X:ballistic}(a)
\label{fig:ballistic.184}}}}
\end{table}

(b) (ECA\#54) We continue
Example \ref{X:variable.width.defect}(b).  In this case, 
$\bL=\bR=\bB$, and $\sL=\sR=\sB:=\{
\white\white\white\black,$ $ \  
\white\white\black\white,$ $ \  
\white\black\white\white,$ $ \  
\black\white\white\white; \ 
\black\black\black\white,$ $ \ 
\black\black\white\black,$ $ \ 
\black\white\black\black,$ $ \ 
\white\black\black\black\}$.  The maps $\varphi_L=\varphi_R$ and
$\varsigma_L=\varsigma_R$ are defined
{\scriptsize
\[
\begin{array}{l}
\white\white\white\black \ \stackrel{\varsigma}{\mapsto} \ 
\white\white\black\white \ \stackrel{\varsigma}{\mapsto} \ 
\white\black\white\white \ \stackrel{\varsigma}{\mapsto} \  
\black\white\white\white \ \stackrel{\varsigma}{\mapsto} \  
\white\white\white\black
\\
\black\black\black\white \ \stackrel{\varsigma}{\mapsto} \  
\black\black\white\black \ \stackrel{\varsigma}{\mapsto} \  
\black\white\black\black \ \stackrel{\varsigma}{\mapsto} \  
\white\black\black\black \ \stackrel{\varsigma}{\mapsto} \  \black\black\black\white\\
\white\white\white\black \ \stackrel{\varphi}{\mapsto} \  
\black\white\black\black \ \stackrel{\varphi}{\mapsto} \  
\white\black\white\white \ \stackrel{\varphi}{\mapsto} \ 
\black\black\black\white \ \stackrel{\varphi}{\mapsto} \  
\white\white\white\black  
\\
\white\white\black\white \ \stackrel{\varphi}{\mapsto} \  
\white\black\black\black \ \stackrel{\varphi}{\mapsto} \  
\black\white\white\white \ \stackrel{\varphi}{\mapsto} \ 
\black\black\white\black \ \stackrel{\varphi}{\mapsto} \  
\white\white\black\white 
\end{array}
\]
}
Consider the $\gam^\pm$ defects in Figure \ref{fig:54.dislocate}.
In this case, $\sD=\sA=\{0,1\}$, so $\sX=\sB\x\sA\x\sB$.
The $\gam^\pm$ particle types correspond to 2-periodic orbit classes 
$\Gam^+$ and $\Gam^-$ of the 
dynamical system $\xi:\sX\into\sX$, where
\beq
 \Gam^+ & := & \lb\{ 
(\mbox{\scriptsize$\white\white\white\black, \,\white,  \,\black\black\black\white$}) \ , \ 
 (\mbox{\scriptsize$\white\black\black\black, \,\white, \,\white\white\black\white$}) \rb\} 
\\
\And \Gam^- & := & \lb\{  
(\mbox{\scriptsize$\white\black\white\white, \, \white, \,\black\white\white\white$ } ) 
\ , \ 
(\mbox{\scriptsize$\white\black\black\black, \, \white, \,\black\black\black\white$  }) \rb\}.
\eeq
Note that $\bd=\zd_0$ is always $\white$. 
Also, $\vV\equiv+1$ on $\Gam^+$, 
so that $\xi(\zl,\white,\zr)= 
(\varsigma\circ\varphi(\zl),\white,\varsigma\circ\varphi(\zr))$
for both $(\zl,\white,\zr)\in\Gam^+$.
Likewise, $\vV\equiv-1$ on $\Gam^-$, 
so that $\xi(\zl,\white,\zr)= 
(\varsigma^{-1}\circ\varphi(\zl),\white,\varsigma^{-1}\circ\varphi(\zr))$
for both $(\zl,\white,\zr)\in\Gam^-$.
}

\Remark{Theorem \ref{thm:ballistic} can easily be generalized to
defect particles in $\ZD$, where $\bX\subset\AZD$ is a transitive,
$\shift{}$-periodic subshift of finite type.  However, if $D\geq 2$,
then such particles cannot be essential defects [see Remark
\ref{defect.particle.remarks}(d)], because if $\bX$ is
$\shift{}$-periodic and $D\geq2 $, then any finite defect in
$\bX$ is removable. Defect particles may still be $\Phi$-persistent,
however.  The most familiar examples of removable, yet persistent,
ballistic defect particles are the `gliders' and `oscillators' of
Conway's {\em Game of Life} \cite{Eppstein,Gotts,BerlekampConwayGuy}
 and its variants
\cite{Bays6,Bays5,Bays4,Bays3,Bays2,Bays1,Evans1,Evans2,Evans3,Evans4,Evans5}.}

\section{The Diffusive Regime \label{S:diffuse}}
  Under certain conditions, a defect particle performs a generalized
random walk.  To demonstrate this, we first review some elementary
probability theory.

\parag{Bernoulli Measures and (hidden) Markov Measures:}
 Let $\sA$ be a discrete set (finite or countable), and let $\sM(\AN)$
be the set of Borel probability measures on $\AN$.  
If $\mu\in\sM(\AN)$, then 
$\mu$ is {\dfn $\shift{}$-invariant} if
$\shift{}(\mu)=\mu$, where $\shift{}(\mu)\in\sM(\AN)$ is 
defined by
$\shift{}(\mu)[\bB]:=\mu[\shift{-1}(\bB)]$
 for any Borel subset $\bB\subset\AN$.
The measure-preserving dynamical system $(\AN,\mu,\shift{})$ is then
called a {\dfn stationary stochastic process}.
For any $m,n\in\Natur$ and any $\bc\in\sA^\CC{0..m}$, let $[\bc]_n:=
\set{\ba\in\AN}{\ba_{\CC{n\ldots n+m}}=\bc}$ be the {\dfn cylinder set}
defined by $\bc$ at position $n$.  
Clearly, $\mu$ is $\shift{}$-invariant iff
$\mu([\bc]_n)=\mu([\bc]_0)$ for all $m,n\in\Natur$ and
$\bc\in\sA^\CC{0..m}$.  Thus, we write ``$\mu[\bc]$'' to mean
$\mu([\bc]_0)$.  
We call $\mu$ a {\dfn Bernoulli measure} if
there is a measure $\mu^0\in\sM(\sA)$ (the `one-point marginal' of $\mu$)
such that, for any $n\in\Natur$
for any $\bc\in\sA^\CC{0..n}$, \ $\mu[\bc] \ = \
\mu^0(\zc_0)\mu^0(\zc_1)\cdots\mu^0(\zc_n)$. 
For example, if $\#(\sA)=A$, then
the {\dfn uniform measure} is the Bernoulli measure $\eta$
with $\eta^0(\za)=\frac{1}{A}$ for all $\za\in\sA$; hence
$\eta(\bc)=\frac{1}{A^{n+1}}$ for all $\bc\in\sA^\CC{0..n}$.

 A measure $\mu\in\sM(\AN)$ 
is a {\dfn Markov measure} if there is a measure $\mu^0\in\sM(\sA)$ 
and a {\dfn transition probability
function} $\tau:\sA\into\sM(\sA)$ such that,
for any $\bc\in\sA^\CC{0..m}$, \ $\mu([\bc]_0) \ = \
\mu^0(\zc_0)\tau(\zc_0,\zc_1)\tau(\zc_1,\zc_2) \cdots\tau(\zc_{n-1},\zc_n)$; see
\cite[\S6.2]{Kitchens}, \cite[\S2.3]{LindMarcus} or
\cite[\S4.4]{Varadhan}.  In this case, $\mu$ is $\shift{}$-invariant
iff $\mu^0$ is {\dfn stationary}, meaning that $\mu^0(b) = \sum_{\za\in\sA} 
\mu^0(\za)\cdot \tau(\za,b)$ for all $b\in\sA$.
For example, any Bernoulli measure is a $\shift{}$-invariant Markov measure,
with $\tau(\za,\zb)=\mu^0(\zb)$ for all $\za,\zb\in\sA$.
If $\sB$ is another set, and  $\psi:\sA\into\sB$ is any function,
we define $\psi^\Natur:\AN\into\BN$ by
$\psi^\Natur(\za_1,\za_2,\ldots):=(\psi(\za_1),\psi(\za_2),\ldots)$.  We say
$\nu\in\sM(\BN)$ is a {\dfn hidden Markov measure} if $\nu =
\psi^\Natur(\mu)$, for some Markov measure $\mu\in\sM(\AN)$ and
function $\psi:\sA\into\sB$.  
Bernoulli/Markov measures on $\AZ$ are defined analogously.

\parag{Random Walks:}
  Let $\dV\subset\Zahl$, and let $\nu\in\sM(\dV^\Natur)$ be a
hidden Markov measure.  Define $\Sigma:\dV^\Natur\into\Zahl^\Natur$ by
$\Sigma(\fv_1,\fv_2,\fv_3,\ldots) \ := \ 
 (0,\fv_1, \ \fv_1+\fv_2, \ \fv_1+\fv_2+\fv_3, \ \ldots)$.
 The probability measure $\omg:=\Sigma(\nu)\in\sM[\Zahl^\Natur]$
is called a (generalized) {\dfn random walk}, 
with {\dfn increment process} $\nu$.
For example, the one-dimensional {\dfn Simple Random Walk} (SRW) is
obtained by setting $\dV:=\{-1,1\}$ and
letting $\nu$ be the Bernoulli measure with $\nu[\pm1]=\frac{1}{2}$;
see \cite[Example 4.1]{Varadhan}.

\ignore{We say $\omg$ is {\dfn commensurable} if $\nu$ is commensurable.
For example, the simple random walk is commensurable.}

\parag{Resolving subshifts:}
 Let $\sB\subseteq\sA$, and let $\bS\subset\BZ\subseteq\AZ$ be a
Markov subshift.   For any $\zb\in\sB$, let
$\sP_\bS(\zb):=\set{\za\in\sB}{(\za,\zb)\in\bS_2}$ be the {\dfn
predecessor set} of $\zb$, and let
$\sF_\bS(\zb):=\set{\zc\in\sB}{(\zb,\zc)\in\bS_2}$ be the {\dfn follower
set} of $\zb$.  We say that $\bS$ is {\dfn left-regular} if there is
some constant $P_\bS\in\Natur$ such that $\#\lb[\sP_\bS(\zb)\rb]=P_\bS$ for all
$\zb\in\sB$.  Likewise $\bS$ is {\dfn right-regular}  if there is some
constant $F_\bS\in\Natur$ such that $\#\lb[\sF_\bS(\zb)\rb]=F_\bS$ for all
$\zb\in\sB$. 

 The {\dfn Parry measure} $\eta\in\sM(\bS)$ is the measure of maximal
$\shift{}$-entropy on $\bS$, and is a Markov measure on $\bS$ which assigns roughly equal
probability to all $\bS$-admissible paths of any given length;
see \cite[Thm.10]{Parry}, \cite[\S13.3]{LindMarcus},
 or \cite[Thm.6.2.20]{Kitchens}.
If $\bS$ is left- or right-regular, then 
$\eta^0$ is the uniform measure on $\sB$.
If $\bS$ is right-regular, then $\tau(\zb,\bullet)$ is the uniform measure on
$\sF(\zb)$ for every $\zb\in\sB$;
that is, $\tau(\zb,\zc)=1/F_\bS$ for all $\zc\in\sF_\bS(\zb)$.
Likewise, if $\bS$ is left-regular, and we define 
the `backwards' transition probability $\backtau:\sA\into\sM(\sA)$ by
$\backtau(\za,\zb):=\eta[\za\zb]/\eta^0(\zb)$, then $\backtau(\bullet,\zb)$ is the
uniform measure on $\sP_\bS(\zb)$ for every $\zb\in\sB$; that is
$\backtau(\za,\zb)=1/P_\bS$ for all $\za\in\sP_\bS(\zb)$.

Let $\Phi:\AZ\into\AZ$ be a CA with
$\Phi(\bS)\subseteq\bS$.  Suppose $\Phi$ has local rule
$\phi:\sA^{\{-1,0,1\}}\into\sA$.  Then $\bS$ is a {\dfn left-resolving
subshift} for $\Phi$ if, for any fixed $(\zb,\zc,\zd)\in\bS_{3}$,
with $\ze:=\phi(\zb,\zc,\zd)$, the
function $\sP_\bS(\zb)\ni \za \mapsto \phi(\za,\zb,\zc)\in \sP_\bS(\ze)$ is injective
\cite[Defn.8.1.7]{LindMarcus}.
If $\bS$ is  left-regular, then `injective' implies `bijective'.
In this case, for any $(\zb,\zc)\in\bS_{2}$, define
$\phi(\bS,\zb,\zc):=\set{\phi(\za,\zb,\zc)}{\za\in\sP_\bS(\zb)}$; 
 then $\#\lb[\phi(\bS,\zb,\zc)\rb]=P_\bS$.
Likewise, $\bS$ is a {\dfn right-resolving subshift} for $\Phi$ if, for
any fixed $(\za,\zb,\zc)\in\bS_{3}$ with $\ze:=\phi(\za,\zb,\zc)$,
the function $\sF_\bS(\zc)\ni \zd \mapsto
\phi(\zb,\zc,\zd)\in \sF_\bS(\ze)$ is injective. 
If $\bS$ is right-regular, then `injective' implies
`bijective'.  In this case, for any $(\zb,\zc)\in\bS_{2}$, define
$\phi(\zb,\zc,\bS):=\set{\phi(\zb,\zc,\zd)}{\zd\in\sF(\zc)}$; then 
$\#\lb[\phi(\zb,\zc,\bS)\rb]=F_\bS$.
If $\bS$ is either left- or right-resolving for $\Phi$, then
$\Phi(\bS)=\bS$, which implies that
the Parry measure on $\bS$ is $\Phi$-invariant \cite[Thm.6.2.21]{Kitchens}.

\example{\label{X:resolve}
(a) Let $\sB\subseteq\sA$ and let $\bS:=\BZ$. 
Then $\BZ$ is left- and right-regular, because $\sP_\bS(\zb)=\sB=\sF_\bS(\zb)$ for
all $\zb\in\sB$. 
If $\Phi:\AZ\into\AZ$ and $\Phi(\BZ)\subseteq\BZ$, then
$\BZ$ is left-resolving for $\Phi$ iff 
$\Phi$ is {\dfn left-permutative} on $\sB$, i.e.
for any  $(\zb,\zc)\in\sB^2$, the
function $\sB \ni \za \mapsto \phi(\za,\zb,\zc)\in \sB$ is bijective \cite{Hedlund}.
Likewise, $\BZ$ is right-resolving iff
$\Phi$ is {\dfn right-permutative} on $\sB$, i.e.
for any  $(\za,\zb)\in\sB^2$, the
function $\sB \ni \zc \mapsto \phi(\za,\zb,\zc)\in \sB$ is bijective.
 The Parry measure is the uniform 
measure on $\BZ$, and is preserved by any permutative cellular automaton.
In the terminology of \cite{Elo93a,Elo93b}, $\sB$ is called a
{\dfn permutive subalphabet} for $\Phi$.

(b) For example, let $(\sB,+)$ be a finite abelian group and let $\bS:=\BZ$.
Then $\Phi$ is a
{\dfn linear cellular automaton} on $\BZ$ if there are endomorphisms
$\varphi_{-1},\varphi_0,\varphi_1\in\End{\sB}$ such that, for all
$\zb_{-1},\zb_0,\zb_1\in\sB$, we have $\phi(\zb_{-1},\zb_0,\zb_1) =
\varphi_{-1}(\zb_{-1}) + \varphi_0(\zb_0)+ \varphi_{1} (\zb_{1})$.  In this
case, $\Phi$ is left- (resp. right-) permutative on $\sB$ iff
$\phi_{-1}$ (resp. $\phi_1$) is an automorphism of $\sB$.  (Note that
we do not require that $\sA$ be a group, or that $\Phi$ be linear on
the rest of $\AZ$.)  Under pointwise addition,
$\BZ$ is a compact abelian group, and
the Parry measure (the uniform measure) is the Haar measure on $\BZ$.

(c)   In particular, if $\sB=\Zahlmod{n}$ for some $n\in\Natur$,
then $\Phi$ is linear if there are 
constants $\varphi_{-1},\varphi_0,\varphi_1\in\Zahl$ such that
 $\phi(\zb_{-1},\zb_0,\zb_1) = (\varphi_{-1} \zb_{-1}
+  \varphi_0 \zb_0+ \varphi_{1} \zb_{1})$ mod $n$.
In this case, $\Phi$ is left- (resp. right-) permutative
on $\sB$ iff $\phi_{-1}$ (resp. $\phi_1$) is relatively prime to $n$.
The Haar measure $\eta$ on $\BZ$ is the `natural' invariant measure
for such linear CA.  For example,
if $\mu\in\sM(\BZ)$ is any measure satisfying broad
conditions (e.g. an $N$-step Markov measure with full support), then
$\Phi$ {\dfn asymptotically randomizes} $\mu$, meaning that
$\D \lim_{N\goto\oo} \frac{1}{N}\sum_{n=1}^N
\Phi^n(\mu)  =  \eta$ in the weak* topology on $\sM(\BZ)$;  see
\cite{Lin84,MaassMartinez1,MaassMartinez2,PivatoYassawi1,PivatoYassawi2,PivatoYassawi3}.
Furthermore, if $n$ is prime, then
$\eta$ is the only $\Phi$-invariant, $\shift{}$-ergodic measure with
positive entropy \cite[Thm. 12]{HostMaassMartinez};
for some generalizations and related results, see also \cite{PivatoPermutative,Sablik}.

(d) If $(\sB,+)$ is a finite group, then 
a {\dfn Markov subgroup} is a Markov subshift $\bS\subset\BZ$
which is also a subgroup of $(\BZ,+)$; see \cite{KitchensXDZDG,KitSch2,KitSch1}
and \cite[\S6.3]{Kitchens}.  It follows that
$\sF_\bS(0)$ and $\sP_\bS(0)$ are subgroups of $\sB$ 
(see \cite[Prop.3(ii)]{KitchensXDZDG} or \cite[Lem.6.3.4(ii,iii)]{Kitchens}). 
Furthermore, $\bS$ is left- and right- regular, because 
for any $\zb\in\sB$, $\sF_\bS(\zb)$ is a coset of $\sF_\bS(0)$, and    
$\sP_\bS(\zb)$ is a coset of $\sP_\bS(0)$ 
(see \cite[Prop.3(iii)]{KitchensXDZDG} or \cite[Lem.6.3.4(iv)]{Kitchens}).
The Parry measure of $\bS$ is then the Haar measure on $\bS$ as a compact group.

If $\Phi$ is a linear CA with local rule $\phi(\zb_{-1},\zb_0,\zb_1) = \varphi_{-1} \zb_{-1}
+  \varphi_0 \zb_0+ \varphi_{1} \zb_{1}$ for some 
constants $\varphi_{-1},\varphi_0,\varphi_1\in\Zahl$ [as in Example (c)]
then $\Phi(\bS)\subseteq\bS$.  
Also, $\Phi$ is left- (resp. right-) permutative
on $\sB$ iff $\phi_{-1}$ (resp. $\phi_1$) acts injectively on
$\sP_\bS(0)$ (resp. $\sF_\bS(0)$).
As in Example (c), many measures on $\bS$ (e.g. Markov measures with
full support) are asymptotically randomized to $\eta$ by such CA; see
\cite{MaassMartinezPivatoYassawi1,MaassMartinezPivatoYassawi2}.

(e) Let $\sB=\{\zb\}$ and let $\bS=\BZ =\{\bb\}$,
 where $\bb=[\ldots \zb \zb \zb \ldots]$. Then $\bS$ is trivially left- and right-regular, because $\sF_\bS(\zb)=\{\zb\}=\sP_\bS(\zb)$.   If $\Phi(\bb)=\bb$, then $\bS$ is 
left- and right- resolving for $\Phi$.  The Parry measure of $\bS$ is the point mass
on $\bb$.
}
Let $\Phi:\AZ\into\AZ$ be a cellular automaton.  A {\dfn resolving system}
for $\Phi$ is a quadruple $(\bL,\bR;\lam,\rho)$, where:
\blist
  \item $\bL,\bR\subset\AZ$ are Markov subshifts, and 
$\bL\union\bR$ is also a Markov subshift.
  \item  $\bL$ is left-regular, $\Phi(\bL)\subseteq\bL$, and $\bL$
is left-resolving for $\Phi$.
  \item $\bR$ is right-regular, $\Phi(\bR)\subseteq\bR$, and $\bR$
is right-resolving for $\Phi$.
  \item $\lam\in\sM(\bL)$ is the Parry measure on $\bL$, and
 $\rho\in\sM(\bR)$ is the Parry measure on $\bR$.
\elist
\example{\label{X:res.sys}If $\sL,\sR\subset\sA$, then $\bA:=\LZ\union\RZ$ is a
subshift of finite type iff either $\sL=\sR$ or they are disjoint.
In this case,  $(\LZ,\RZ;\lam,\rho)$ is a resolving system
for $\Phi$ iff:
[i] \ $\lam$ (resp. $\rho$) is the uniform measure on
$\LZ$ (resp. $\RZ$);  \ [ii] \  
$\Phi(\LZ)=\LZ$ and $\Phi$ is left-permutative on $\sL$; and \ [iii] \
$\Phi(\RZ)=\RZ$ and $\Phi$ is right-permutative on $\sR$,
as in Example \ref{X:resolve}(a).}

Let $(\bL,\bR;\lam,\rho)$ be a resolving system, let $\bA:=\bL\union\bR$, and
let $\bD^{W,0}_{\bL,\bR}$ be
the set of all elements in $\AZ$ with a single 
$(\bL,\bR)$-defect of width $W$ at zero.  Let $\sD:=\sA^W$.
If $\del\in\sM(\sD)$, then we regard $\lam\tensor\del\tensor\rho$ as a
probability measure on $\bD^{W,0}_{\bL,\bR}$ in the obvious way.
Define $\zeta:\bD^W_{\bL,\bR}\into\Zahl^\Natur$
by $\zeta(\ba):=(z_0,z_1,z_2,\ldots)$, where, for all $t\in\Natur$,
\ $z_t=z_t(\ba)\in\Zahl$ is as in 
eqn.(\ref{variable.width.defect.2}). In other words, $\zeta(\ba)$
tracks the trajectory of the defect particle over time.
If $\mu\in\sM(\bD^{W,0}_{\bL,\bR})$,
then $\zeta(\mu)$ is a probability measure on $\Zahl^\Natur$.
The main result of this section is:

\Theorem{\label{defect.random.walk}}
{
  Let $\Phi:\AZ\into\AZ$ be a CA, and let
$(\bL,\bR;\lam,\rho)$ be a resolving system for $\Phi$.
Let $W\in\Natur$, let $\del$ be any probability measure on
$\sD:=\sA^W$, and let $\mu:=\lam\tensor\del\tensor\rho
\in\sM(\bD^{W,0}_{\bL,\bR})$.
Then $\omg:=\zeta(\mu)\in\sM(\Zahl^\Natur)$ is a random walk.
}

\begin{figure}
\centerline{
\begin{tabular}{cc}
\includegraphics[height=14em,width=17em]{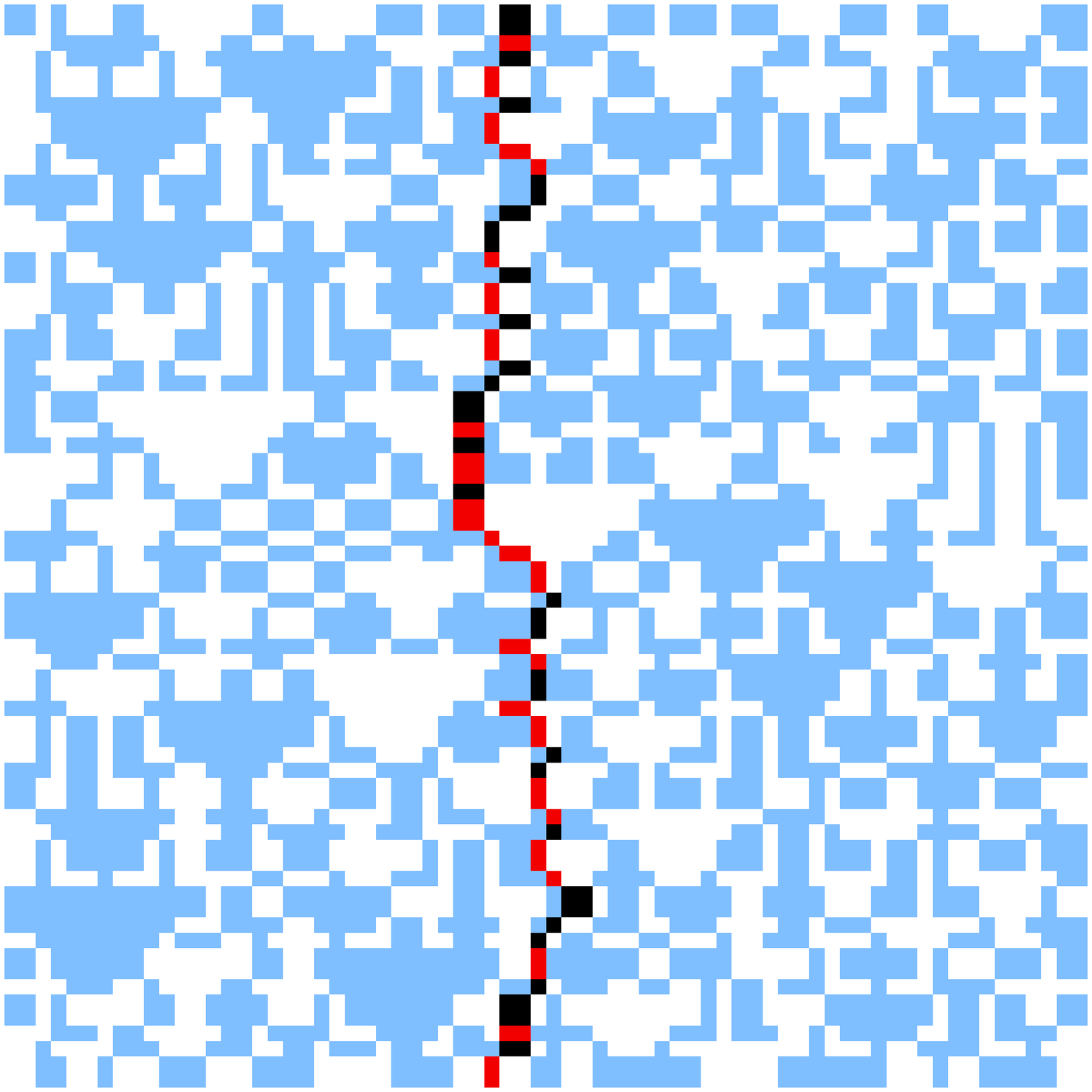} &
\includegraphics[height=14em,width=17em]{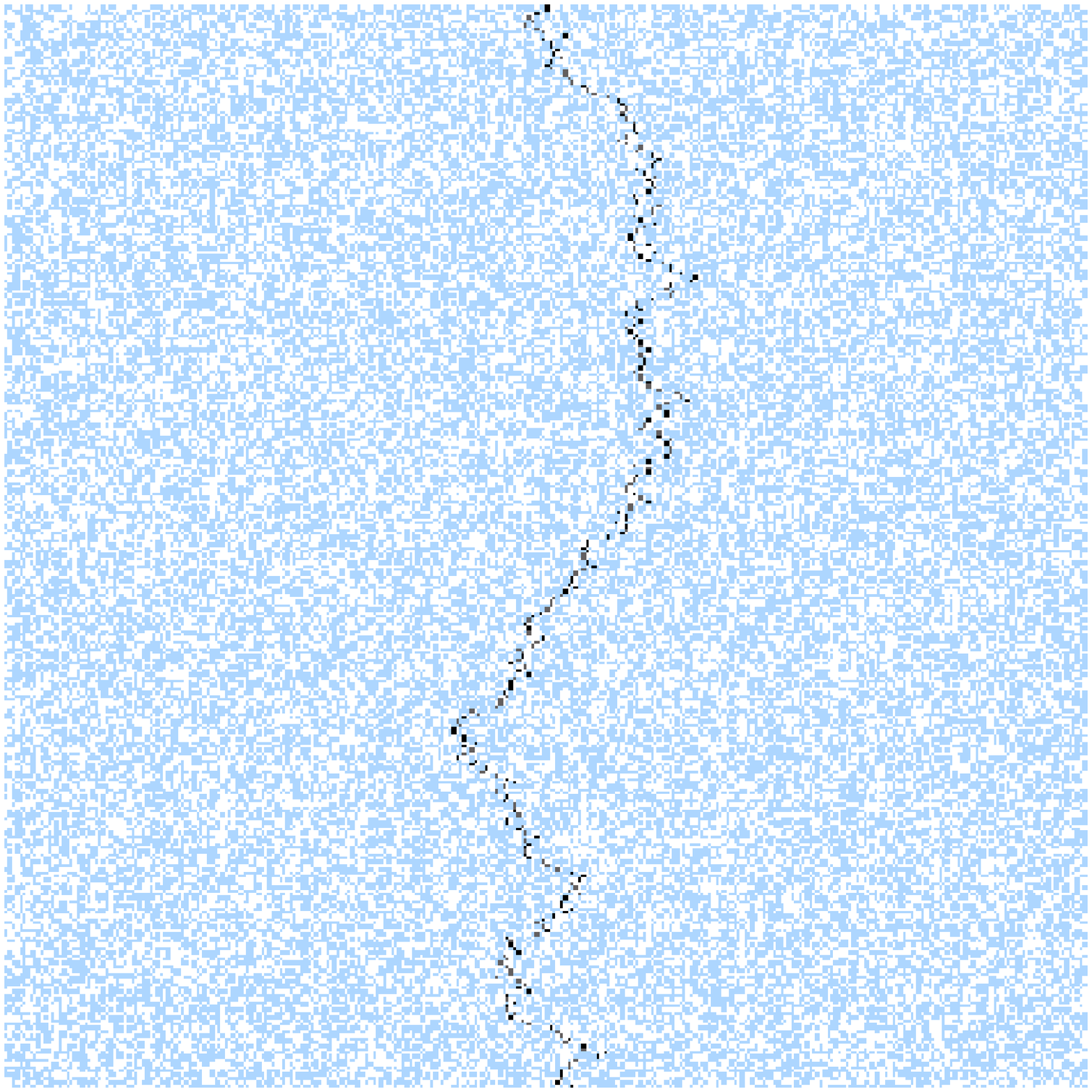} \\
{\scriptsize {\bf(A)}\qquad\qquad Scale: $50 \x 50$} & {\scriptsize 
{\bf(B)}\qquad\qquad Scale: $300 \x 6000$ (space $\x$ time)}\\
\end{tabular}}
\caption{{\footnotesize The randomly walking defect particle of Example 
\ref{X:diffuse}. \label{fig:diffuse3}}}
\end{figure}

%Before proving Theorem \ref{defect.random.walk}, we illustrate with an example.

\example{\label{X:diffuse} Let $\sA:=\Zahlmod{2}\x\{\circ,\bullet\}$.
Define $\phi:\sA^{\{-1,0,1\}}\into\sA$ by
$\phi\lb( {\za_{-1}\atop \zb_{-1}} \ {\za_{0}\atop \zb_{0}} \  {\za_{1}\atop \zb_{1}}\rb)
\ := \ \lb({\za \atop \zb}\rb)
$, where 
{\footnotesize
\beq
\za & :=& \choice{ \za_{-1} + \za_0 + \za_1 &&\If \zb_{-1}=\zb_0=\zb_1= \circ; \\ 
 \za_{-1} + \za_0  &&\If \zb_{-1}=\zb_0=\circ \And \zb_1=\bullet; \\ 
 \za_{0} + \za_1  &&\If \zb_{-1}=\bullet \And \zb_{0}=\zb_1=\circ; \\
 1- \za_0  &&\If \zb_0=\bullet. } 
\\
\zb &:=& \choice{ \bullet &&\If \zb_{-1}=\bullet \And \za_{-1}=\za_0=0; \\
 \bullet && \If  \zb_{1}=\bullet \And \za_{0}=\za_1=1; \\
 \circ && \If  \zb_{0}=\bullet \And \za_0=\za_1=0;\\ 
 \circ && \If  \zb_{0}=\bullet \And \za_{-1}=\za_0=1; \\
 \zb_0 & & \mbox{otherwise.} \\}
\eeq
}
Let $\sL=\sR=\Zahlmod{2}\x\{\circ\}$, which we identify
with $\Zahlmod{2}$;  then $\Phi$ acts on $\bL:=\sL^\Zahl=\sR^\Zahl=:\bR$ like the
linear cellular automaton $\Psi:(\Zahlmod{2})^\Zahl\into(\Zahlmod{2})^\Zahl$
with local rule $\psi(\zx_{-1},\zx_0,\zx_1):=\zx_{-1}+\zx_0+\zx_1 \pmod{2}$.  Thus,
$\sL=\sR$ is a left- and right-permutative subalphabet for $\Phi$
[see Example \ref{X:resolve}(c)]; if $\lam=\rho$ is the 
uniform measure on $\bL=\bR$, then $(\bL,\bR;\lam,\rho)$ is
a resolving system, as in Example \ref{X:res.sys}.
The set $\Zahlmod{2}\x\{\bullet\}$ is the set of defect states.
An element of $\bD^1_{\bL,\bR}$ has the form
$\lb[ \ldots {\zl_3 \atop \circ} \ 
{\zl_2 \atop \circ} \ 
{\zl_1 \atop \circ} \ 
{\zd_0 \atop \bullet} \
{\zr_1 \atop \circ} \ {\zr_2 \atop \circ} \ {\zr_3 \atop \circ} \ldots\rb]$,
where $\zl_i,\zr_i,\zd_0\in\Zahlmod{2}$. 
The defect particle `$\bullet$' moves left if $\zl_1=\zd_0=1$, and moves
right if $\zd_0=\zr_1=0$; otherwise it remains stationary.
Figure \ref{fig:diffuse3}(A) shows a close-up  spacetime
diagram of the resulting random walk, while Figure
\ref{fig:diffuse3}(B) shows a large-scale  spacetime
diagram of the same walk.  
}

If $\mu\in\sM(\AN)$, then we write, ``For $\forall_\mu \ \ba\in\AN,$ \
[statement]'', or ``[statement], \muae'', to mean
``$\mu\set{\ba\in\AN}{\mbox{[statement] is true}}=1$''.  If
$\dI\subset\Natur$, let $\pr{\dI}:\AN\into\sA^\dI$ be the projection
map.  Thus, if $\mu\in\sM(\AN)$, then $\pr{\dI}(\mu)\in\sM(\sA^\dI)$.
If $\dJ\subset\Natur$ and $\bb\in\sA^\dJ$, then let $\mu_{|\bb}$ be
the conditional probability measure on $\AN$ given $\bb$; in other words, for
any $\bU\subset\AN$, $\mu_{|\bb}(\bU) :=
\mu(\bU\intsct[\bb])/\mu[\bb]$, where $[\bb]$ is the cylinder set
defined by $\bb$.  More generally, if $\bU\subset\AN$ is
Borel-measurable, and if $\SigAlg$ is a $\sigma$-subalgebra of the
Borel $\sigma$-algebra on $\AN$, then let $\mu_{|\SigAlg}(\bU)$ be the
conditional probability function of $\bU$ given $\SigAlg$;
i.e. $\mu_{|\SigAlg}(\bU)$ is the $\SigAlg$-measurable function such
that, for any $\bS\in\SigAlg$, $\int_{\bS} \mu_{|\SigAlg}(\bU) \ d\mu
\ = \ \mu[\bS\intsct\bU]$. This function is uniquely defined 
\muae; see e.g. \cite[\S4.3]{Varadhan}.

 In particular, if $\dJ\subset\Natur$, let $\SigAlg(\dJ)$ be
the sigma-algebra generated by all cylinder sets $[\zc]_j$,
where $\zc\in\sA$ and $j\in\dJ$ (hence $\SigAlg(\Natur)$ is the Borel sigma-algebra
of $\AN$). 
A $\shift{}$-invariant $\mu\in\sM(\AN)$ is {Bernoulli} iff, for
any disjoint subsets $\dI,\dJ\subset\Natur$, and any $\bb\in\sA^\dI$, \ 
$\mu_{|\SigAlg(\dJ)}(\bb) \equiv \mu(\bb)$ \muae.
We say $\mu$ is {\dfn Markovian} iff for
any $m\in\Natur$ and $\zb\in\sA$,
 \  $\mu_{|\SigAlg\CC{0\ldots m}}([\zb]_{m+1}) = \mu_{|\SigAlg\{m\}}([\zb]_{m+1})$.
Thus, $\mu$ is a Markov measure if $\mu$ is Markovian and if, furthermore,
for any $\za\in\sA$ and for $\forall_\mu \ \bx\in[\za]_m$, we have
$\mu_{|\SigAlg\{m\}}([\zb]_{m+1})(\bx) \equiv \tau(\za,\zb)$.

\begin{figure}
\centerline{
\psfrag{A}[][]{$\SigAlg\{t\}$}
\psfrag{B}[][]{$\SigAlg\CC{0...t}$}
\psfrag{C}[][]{$\SigAlg^*$}
\psfrag{D}[][]{{\footnotesize $\lam$-random initial condition in $\bL^-$}}
\psfrag{E}[][]{{\footnotesize Incoming $\lam$-random `noise'}}
\psfrag{F}[][]{{\footnotesize $\rho$-random initial condition in $\bR^+$}}
\psfrag{G}[][]{{\footnotesize Incoming $\rho$-random `noise'}}
\includegraphics[height=35em,width=20em,angle=-90]{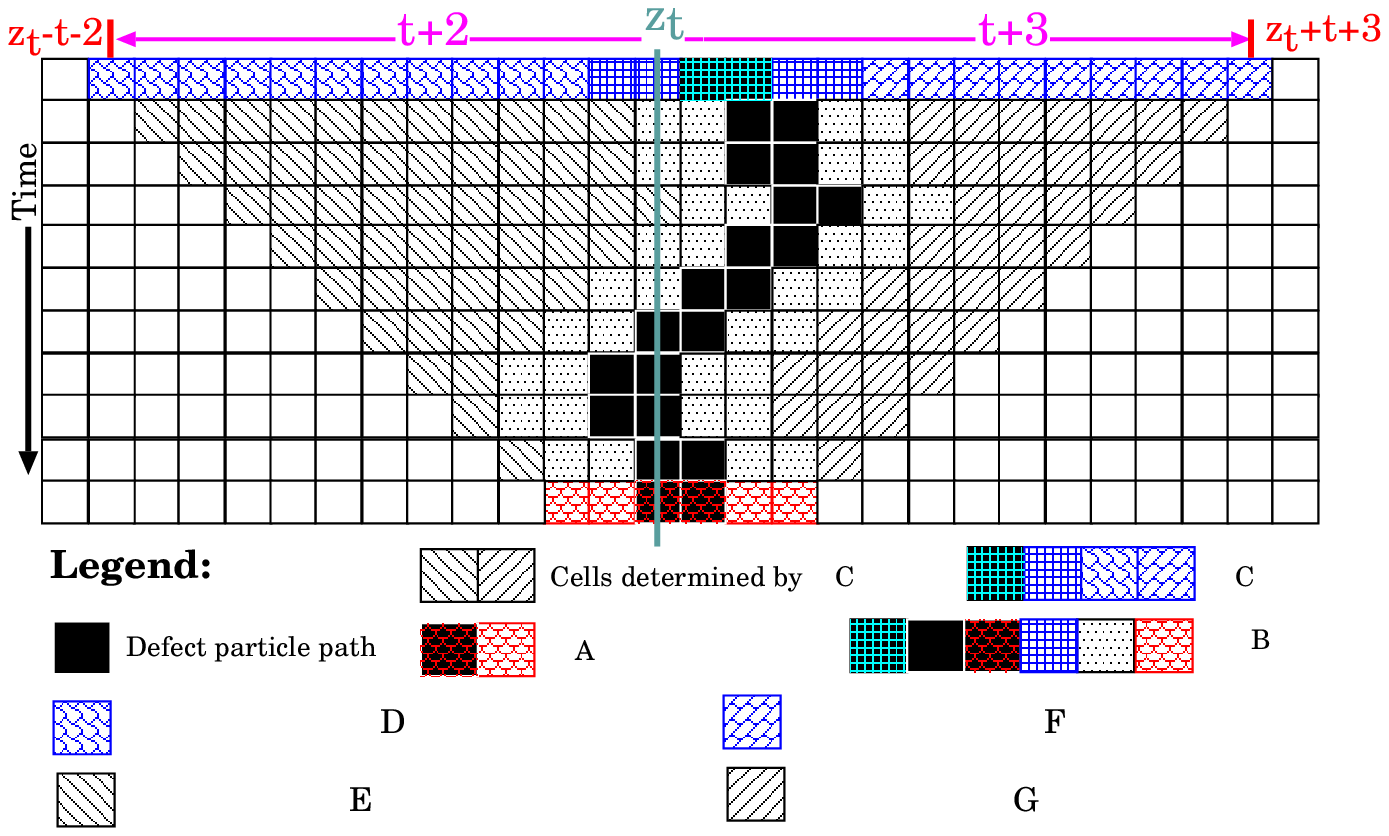}}
\caption{{\footnotesize A schematic spacetime diagram illustrating
the sigma algebras $\SigAlg\{t\}
\subseteq \SigAlg\CC{0...t}\subseteq \SigAlg^*$ in
Theorem \ref{defect.random.walk}. \label{fig:diffuse}}}
\end{figure}

\begin{figure}
\centerline{
\psfrag{a}[][]{{\scriptsize $\boldsymbol{\zl_{2}^t}$}}
\psfrag{b}[][]{{\scriptsize $\boldsymbol{\zl_{1}^t}$}}
\psfrag{c}[][]{{\scriptsize \textcolor{white}{$\boldsymbol{ \zd^t_0}$}}}
\psfrag{d}[][]{{\scriptsize \textcolor{white}{$\boldsymbol{ \zd^t_1}$}}}
\psfrag{e}[][]{{\scriptsize $\boldsymbol{ \zr_{1}^t}$}}
\psfrag{f}[][]{{\scriptsize $\boldsymbol{ \zr_{2}^t}$}}
\psfrag{g}[][]{{\scriptsize {$\boldsymbol{ \zl_{4}^t}$}}}
\psfrag{h}[][]{{\scriptsize {$\boldsymbol{ \zl_{3}^t}$}}}
\psfrag{i}[][]{{\scriptsize $\boldsymbol{ \zl_{2}^{t\!+\!1}}$}}
\psfrag{j}[][]{{\scriptsize $\boldsymbol{ \zl_{1}^{t\!+\!1}}$}}
\psfrag{k}[][]{{\scriptsize \textcolor{white}{$\boldsymbol{ \zd^{t\!+\!1}_0}$}}}
\psfrag{l}[][]{{\scriptsize \textcolor{white}{$\boldsymbol{ \zd^{t\!+\!1}_1}$}}}
\psfrag{m}[][]{{\scriptsize $\boldsymbol{ \zr_{1}^{t\!+\!1}}$}}
\psfrag{n}[][]{{\scriptsize $\boldsymbol{ \zr_{2}^{t\!+\!1}}$}}
\psfrag{o}[][]{{\scriptsize {$\boldsymbol{ \zr_{3}^t}$}}}
\psfrag{p}[][]{{\scriptsize {$\boldsymbol{ \zr_{4}^t}$}}}
\psfrag{G}[][]{{\scriptsize {$\boldsymbol{ \zl_{z\!-\!t\!-\!4}^0}$}}}
\psfrag{H}[][]{{\scriptsize {$\boldsymbol{ \zl_{z\!-\!t\!-\!3}^0}$}}}
\psfrag{O}[][]{{\scriptsize {$\boldsymbol{ \zr_{z\!+\!t\!+\!4}^0}$}}}
\psfrag{P}[][]{{\scriptsize {$\boldsymbol{ \zr_{z\!+\!t\!+\!5}^0}$}}}
\psfrag{A}[][]{$\xi_t$ }
\psfrag{B}[][]{Cells determined by $\bb$}
\psfrag{C}[][]{{\footnotesize Incoming `noise':  cells influenced by
left/right neighbours of $\bb$.}}
\psfrag{D}[][]{$\xi_{t+1}$}
\psfrag{E}[][]{$\bb$}
\psfrag{F}[][]{Defect particle's path}
\includegraphics[height=37em, width=37em]{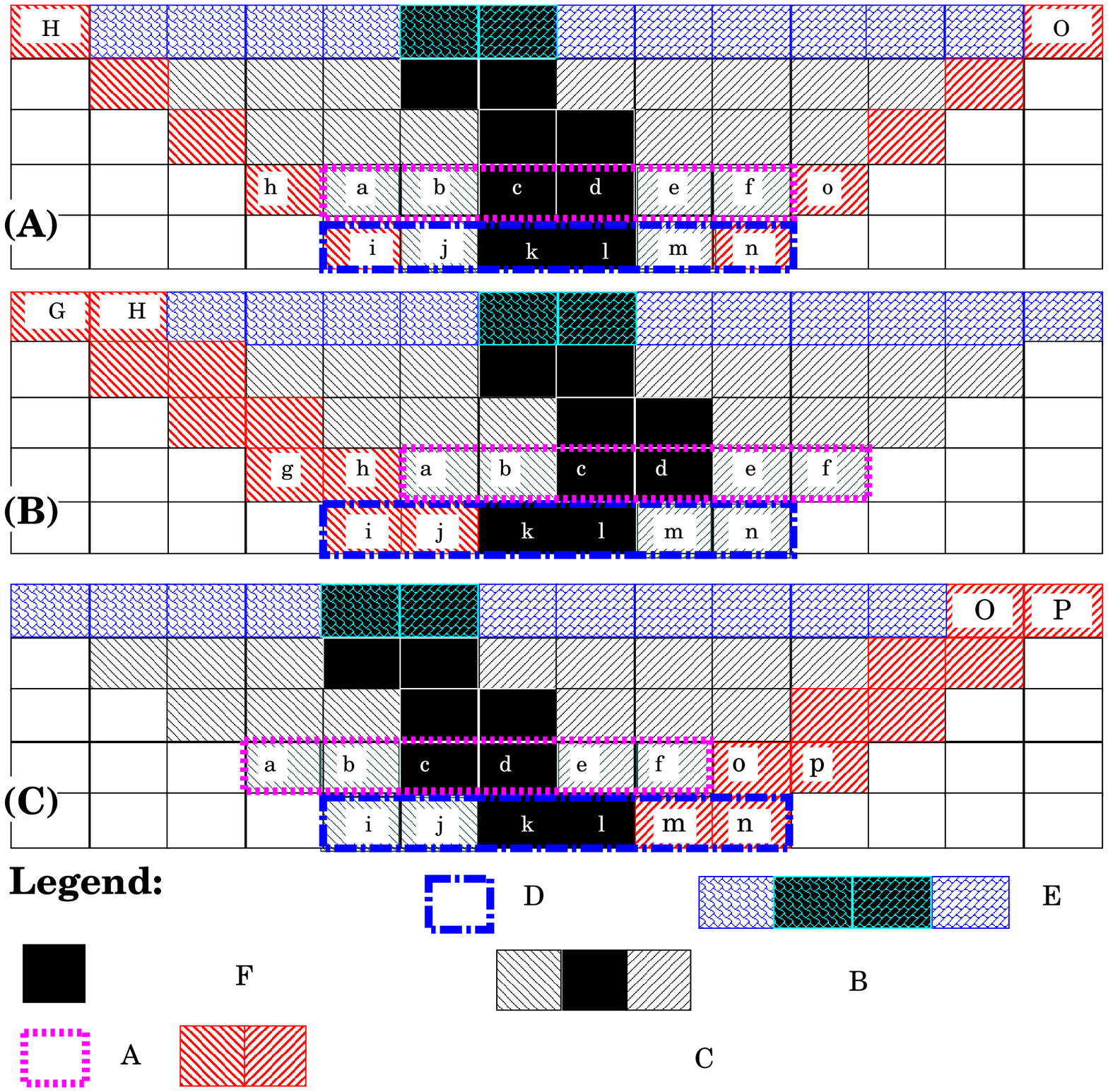}}
\caption{{\footnotesize Claims 1.1 to 1.3 of Theorem \ref{defect.random.walk}. \label{fig:diffuse2}}}
\end{figure}

{\em Proof idea for Theorem \ref{defect.random.walk}:} The left-hand
measure $\lam$ and right-hand measure $\rho$ provide a continual
influx of `random noise'.  The `$\lam$-noise' propagates rightwards
with unit speed because $\bL$ is left-resolving for $\Phi$, whereas
the `$\rho$-noise' propagates leftwards with unit speed because $\bR$
is right-resolving for $\Phi$.  As shown in Figure \ref{fig:diffuse},
the defect particle's trajectory from time $0$ to time $t$ is entirely
determined by the information contained inside of a backwards
`lightcone' emanating from its position at time $t$ back to the
initial state at time zero.  If the particle steps to the left
[respectively, right] at time $t$, then it must step into the path of
incoming $\lam$-noise [respectively, $\rho$-noise] which is outside of
this lightcone, and hence, statistically independent of the particle's
previous trajectory; see Figure \ref{fig:diffuse2}(B)
[respectively, Figure \ref{fig:diffuse2}(C)].  If the particle
stays put at time $t$, then it is exposed to both fresh $\lam$-noise
{\em and} fresh $\rho$-noise; see Figure \ref{fig:diffuse2}(A).
In all three cases, the particle is subjected to fresh perturbations
at time $t+1$ which are statistically independent of its previous
behaviour.  Furthermore, $\lam$ and $\rho$ are $\Phi$-invariant, so
the probability distribution of these perturbations is constant over
time; hence they can be treated as a stationary Markov process, which
drives the particle's motion.

\bthmprf[Proof of Theorem \ref{defect.random.walk}:]
Let $\hsA:=\sA^W$,  let $\hPhi:\hsA^\Zahl\into\hsA^\Zahl$ be the
$W$th-power representation of $\Phi$;  and let $\hgL,\hgR\subset\hsA^\Zahl$
be the $W$th-power representations of $\bL$ and $\bR$.
Then  $\hgL$ (resp. $\hgR$) is still left-regular (resp. right-regular) and
is still left-resolving (resp. right-resolving) for $\hPhi$.
Thus, we
can replace $\sA$ with $\hsA$, $\bL$ with $\hgL$, and $\bR$ with
$\hgR$ and proceed.
By Remark \ref{defect.particle.remarks}(c),
we can thus assume that $W=2$ and that
$\Ups:\bL_1\x\sA^2\x\bR_1\into\sA^2$ and
$\vV:\bL_1\x\sA^2\x\bR_1\into\{-1,0,1\}$. A generic element of 
$\bD^2_{\bL,\bR}$ has the form
\[
 \ba\quad=\quad
[\ldots \zl_{3} \ \zl_{2} \  \zl_{1} \ \zd_0 \ \zd_1 \ \zr_{1}
 \ \zr_{3} \ \zr_{3} \ldots],
\]
where $\zl_n:=\za_{z-n}\in\bL_1$ and $\zr_n:=\za_{z+n+2}\in\bR_1$ for all $n\in\Natur$, while $\zd_i:=\za_{z+i}\in\sA$ for $i=0,1$, with $z\in\Zahl$ being the location of the defect.

Let $\sX:=\bL_2\x\sA^2\x\bR_2$, and
 define $\xi:\bD^2_{\bL,\bR} \into \sX$ so that, if $\ba$ is as above, then
$\xi(\ba):=(\zl_{2},\zl_{1}; \ \zd_0,\zd_1; \ \zr_{1},\zr_{2})$.
For any $t\in\Natur$, let $\xi_t:=\xi\circ\Phi^t$.
In other words, $\xi_t(\ba):=(\zl^t_{2},\zl^t_{1}; \ \zd^t_0,\zd^t_1; \ \zr^t_{1},\zr^t_{2})$,
where $\Phi^t(\ba) \ = \
[\ldots \zl^t_{3} \ \zl^t_{2} \  \zl^t_{1} \ \zd^t_0 \ \zd^t_1 \ \zr^t_{1}
 \ \zr^t_{2} \ \zr^t_{3} \ldots]$.
Next, define $\Xi:\bD^{2,0}_{\bL,\bR} \into  \XN$
by $\Xi(\ba) := (\xi_0(\ba),\xi_1(\ba),\xi_2(\ba),\ldots)$.
Clearly, $\Xi\circ\Phi=\shift{}\circ\Xi$.
Recall that $\vV$ is a function from $\bL_1\x\sA^2\x\bR_1$ into $\{-1,0,1\}$;
treat this as a function $\vV:\sX\into\{-1,0,1\}$, and
apply it coordinatewise to define $\vV^\Natur:\XN\into\{-1,0,1\}^\Natur$.

If $\barmu:=\Xi(\mu)\in\sM(\XN)$,
and $\nu:=\vV^\Natur(\barmu)\in\sM(\{-1,0,1\}^\Natur)$,
then $\omg = \Sigma(\nu)$.
Hence, if $\barmu$ is Markov, then $\nu$ is hidden Markov, so that
$\omg$ is a random walk, as desired.
It remains to show that $\barmu$ is a Markov measure.

Fix $t\in\Natur$.  Let $\SigAlg\{t\}$ be the sigma-algebra on $\AZ$
generated by $\xi_t$, and let $\SigAlg\CC{0...t}$ 
be the sigma-algebra on $\AZ$ generated by $(\xi_0,\xi_1,\ldots,\xi_t)$.
For any $\zx\in \sX$, let
$\bU^{t+1}_\zx:=\xi_{t+1}^{-1}\{\zx\}$
be the set of all initial conditions in $\bD^{2,0}_{\bL,\bR}$ such
that the defect particle at time $t+1$ has internal state
$\zx$.  To show that $\barmu$ is a Markov measure, we must 
find some transition probability function
$\tau:\sX\into\sM(\sX)$ such that:
\beqn
\label{markov.eqn}
\mbox{For all $\zx\in\sX$ and $t\in\Natur$}, \ \
  \mu_{|\SigAlg\CC{0\ldots t}}\lb[\bU^{t+1}_\zx\rb]
\ = \
  \mu_{|\SigAlg\{t\}}\lb[\bU^{t+1}_\zx\rb]
\ = \ \tau(\xi_t,\zx).
\eeqn
For any $z\in\Zahl$, let 
$\bD^t_z:=\set{\ba\in\bD_{\bL,\bR}^{2,0}}{z_t(\ba)=z}$,
and let $\SigAlg^*_z$ be the sigma-algebra on $\bD^t_z$ generated by cylinder
sets in coordinates $\CC{z\!-\!t\!-\!2 \ldots z\!+\!t\!+\!3}$.
Then let $\SigAlg^*$ be the sigma-algebra on $\bD_{\bL;\bR}^{2,0}$ generated by
$\Union_{z\in\Zahl} \SigAlg^*_z$.
Clearly, $\SigAlg\CC{0..t}\subseteq  \SigAlg^*$, because the information
contained in $\SigAlg^*$ is sufficient to determine the first $t$
positions $(z_1,\ldots,z_t)$ of the defect particle,
and its first $t$ internal states $(\xi_1,\ldots,\xi_t)$;
see Figure \ref{fig:diffuse}.

\Claim{There exists a function  $\tau:\sX\into\sM(\sX)$ such that,
for any $\zx\in\sX$ and $t\in\Natur$, we have
$\mu_{|\SigAlg^*}[\bU^{t+1}_\zx] \ = \ \mu_{|\SigAlg\{t\}}[\bU^{t+1}_\zx]
\ = \ \tau(\xi_t,\zx)$.}
\bclaimprf
Let $\zx:=(\zl_2,\zl_1;\zd_0,\zd_1,\zr_1,\zr_2)$, where
$(\zl_2,\zl_1)\in\bL_2$, $(\zd_0,\zd_1)\in\sA^2$, and
$(\zr_1,\zr_2)\in\bR_2$.  Fix $t\in\Natur$.  Let
$\xi_{t}=(\zl^{t}_2,\zl^{t}_1;$ $\ \zd^{t}_0,\zd^{t}_1;$ $ \
\zr^{t}_0,\zr^{t}_1)$ and $\xi_{t+1}=(\zl^{t+1}_2,\zl^{t+1}_1;$ $ \
\zd^{t+1}_0,\zd^{t+1}_1;$ $ \ \zr^{t+1}_0,\zr^{t+1}_1)$, where we regard these
as twelve measurable functions on $\bD_{\bL;\bR}^{2,0}$.  For
$\fv\in\{-1,0,1\}$, let $\bD^{\fv} := \set{\ba\in\bD_{\bL,\bR}^{2,0}}{
\vV(\xi_t(\ba))=\fv}$, and let $\mu^{\fv}:=\mu\restr{\bD^{\fv}}$.
Then $\bD_{\bL,\bR}^{2,0}= \bD^{-1}\disj\bD^{0}\disj\bD^{1}$, and $\mu
= \sum_{\fv=-1}^1 \mu[\bD^{\fv}]\cdot \mu^{\fv}$.  We will thus
consider $\mu^{(-1)}$, $\mu^{0}$ and $\mu^{1}$ separately.

For $\fv\in\{-1,0,1\}$ and $z\in\Zahl$,
let $\bD^\fv_z:=\set{\ba\in\bD^\fv}{z_t(\ba)=z}$, and let
$\bB_z^\fv:=\pr{\CC{z\!-\!t\!-\!2 \ldots z\!+\!t\!+\!3}}
[\bD^\fv_z] \subset \sA^\CC{z\!-\!t\!-\!2 \ldots z\!+\!t\!+\!3}$.
Then let $\bB^\fv:=\D \Disj_{z\in\Zahl} \bB^\fv_z$.  

\setcounter{subclaimcount}{-1}
\subclaim{For any $\bb\in\bB^\fv_z$, let $[\bb]:=\set{\ba\in\bD_{\bL,\bR}^{2,0}}
{\ba_\CC{z\!-\!t\!-\!2 \ldots z\!+\!t\!+\!3}=\bb}$.
Then $\bD^\fv =\D \Disj_{\bb\in\bB^\fv} [\bb]$.}
\bsubclaimprf
$\bD_z^\fv  =  \D\Disj_{\bb\in\bB_z^\fv} [\bb]$,
for any $z\in\Zahl$.
Thus,
$\D\bD^\fv   = 
\Disj_{z\in\Zahl}  \bD_z^\fv
 = 
\Disj_{z\in\Zahl}  \Disj_{\bb\in\bB_z^\fv} [\bb]
 = 
\Disj_{\bb\in\bB^\fv} [\bb]$.
\ignore{
``$\D \Disj_{\bb\in\bB^\fv} [\bb]\subseteq \bD^\fv$'': \ 
If $\bb\in\bB^\fv_z$, then $[\bb]\subseteq\bD_z^\fv$, because 
$\ba_\CC{z\!-\!t\!-\!2 \ldots z\!+\!t\!+\!3}$ is sufficient to determine the
values of $z_t(\ba)$ and $\xi_t(\ba)$.

``$\bD^\fv \subseteq \D \Disj_{\bb\in\bB^\fv} [\bb]$'': \ 
Suppose $\fv=-1$.  
Note that $\D (\Disj_{\bb \in \bB^{(-\!1)}} [\bb] ) \disj
(\Disj_{\bb \in \bB^{0}} [\bb] ) \disj
(\Disj_{\bb \in \bB^{1}} [\bb] )
 =  \bD_{\bL,\bR}^{2,0}  = 
\bD^{(-\!1)}\disj\bD^{0}\disj\bD^{1}$.  Since
$\D\Disj_{\bb\in\bB^0} [\bb]\subseteq \bD^0$
and $\D\Disj_{\bb\in\bB^1} [\bb]\subseteq \bD^1$,
we must have 
$\bD^{(-\!1)} \subseteq \D \Disj_{\bb\in\bB^{(-\!1)}} [\bb]$.  The same reasoning
applies for $\fv=0$ and $\fv=1$.}
\esubclaimprf

\subclaim{For any $\zx\in\sX$, \ 
$\mu^0_{|\SigAlg^*}[\bU^{t+1}_\zx]  =  \mu^0_{|\SigAlg\{t\}}[\bU^{t+1}_\zx]
 =  \tau_{0}\lb({\xi_t\atop \zx}\rb)$, where

\centerline{$\D
\tau_0\lb({{\zl^t_2,\zl^t_1;\zd^t_0,\zd^t_1;\zr^t_1,\zr^t_2}
\atop {\zl_2,\zl_1;\zd_0,\zd_1;\zr_1,\zr_2}}\rb)
\ \ := \ \ \choice{ \D\frac{1}{P_\bL F_\bR} & \If &
\mbox{$\footnotesize \begin{array}[t]{l} 
\zl_2\in\sP_\bL(\zl_1), \ 
\zl_1 = \phi(\zl^t_2,\zl^t_1,\zd^t_0), \\ 
\zd_0= \phi(\zl_{1}^t,\zd^t_0,\zd^t_1), \ \
\zd_1 = \phi(\zd^t_0,\zd^t_1,\zr^t_1), \\
\zr_1= \phi(\zd^t_1, \zr_{1}^t,\zr_2^2),
\And \zr_2\in\sF_\bR(\zr_2);
\end{array}$} \\
0 &&\mbox{otherwise.}}$}}
\bsubclaimprf
Figure \ref{fig:diffuse2}(A) shows how
the values of $(\zl^{t+1}_1,\zd^{t+1}_0,\zd^t_1; \ \zr_{1}^{t+1})$ are
determined by the data in $\SigAlg\{t\}$, because
$\zl_{1}^{t+1} = \phi(\zl^t_2,\zl^t_1,\zd^t_0)$, \ 
$\zd^{t+1}_0 = \phi(\zl_{1}^t,\zd^t_0,\zd^t_1)$, \
$\zd^{t+1}_1 = \phi(\zd^t_0,\zd^t_1,\zr^t_1)$, 
and $\zr_{1}^{t+1} = \phi(\zd^t_1, \zr_{1}^t,\zr_2^2)$.
However, $\zl_{2}^{t+1} = \phi(\zl_{3}^t,\zl_{2}^t,\zl_{1}^t)$ 
and $\zr^{t+1}_{2} = \phi(\zr_{1}^t,\zr^t_{2},\zr^t_3)$
not determined, even by $\SigAlg^*$.
Instead, for any fixed $z\in\Zahl$ and $\bb\in\bB_z^0$,
there is a function $\Phi_\bb:
\sP_\bL(\zl_{z\!-\!t\!-\!2}^{0}) \x \sF_\bR(\zr_{z\!+\!t\!+\!3}^{0})
\into \sP_\bL(\zl_{1}^{t+1}) \x \sF_\bR(\zr_{1}^{t+1})$
such that 
$(\zl_{2}^{t+1},\zr_{2}^{t+1})
\ = \ \Phi_\bb(\zl_{z-t\!-\!3}^0,\zr_{z\!+\!t\!+\!4}^0)$. 
 Furthermore,
$\Phi_\bb$ is bijective,  because $\bL$ is $\Phi$-left-resolving
and $\bR$ is $\Phi$-right-resolving.

The set $\sP_\bL(\zl_{z\!-\!t\!-\!2}^{0}) \x \sF_\bR(\zr_{z\!+\!t\!+\!3}^{0})$ has cardinality $P_\bL F_\bR$, because  $\bL$ is left-regular and
 $\bR$ is right-regular.
Let $\mu_{|\bb}$ be the conditional measure on $\bD^0$
given $\bb$.
If $\eta^0_\bb:=\pr{\{z\!-\!t\!-\!3,z\!+\!t\!+\!4\}}(\mu_{|\bb})$,
then $\eta^0_\bb$ is
the uniform measure assigning mass $1/(P_\bL F_\bR)$ 
to each element of $\sP_\bL(\zl_{z\!-\!t\!-\!2}^{0}) \x \sF_\bR(\zr_{z\!+\!t\!+\!3}^{0})$, because
$\lam$ is the Parry
measure on $\bL$ and
$\rho$ is the Parry measure on $\bR$.
Note that any $\bb\in\bB^0$ completely determines the values of
$\zl_{1}^{t+1}$ and $\zr_{1}^{t+1}$ (because these
are $\SigAlg^*$-measurable functions).  Let
$\eta^{t+1}_\bb$
be the uniform measure assigning mass $1/(P_\bL F_\bR)$ 
to each element of $\sP_\bL[\zl_{1}^{t+1}(\bb)] \x \sF_\bR[\zr_{1}^{t+1}(\bb)]$.
Then the $\mu^0_{|\bb}$-probability distribution of
 $(\zl_{2}^{t+1},\zr_{2}^{t+1})$
is the measure 
$\Phi_\bb(\pr{\{z\!-\!t\!-\!3,z\!+\!t\!+\!4\}}(\mu^{0}_{|\bb}))
=\Phi_\bb(\eta^0_\bb) \ \eeequals{(*)} \ \eta^{t+1}_\bb$ (here $(*)$ is
because $\Phi_\bb$ is bijective, while 
both $\eta^0_\bb$ and $\eta^{t+1}_\bb$ are uniform measures
on sets with $P_\bL F_\bR$ elements).  Thus,
$\mu^0_{|\bb}[\bU^{t+1}_\zx] =  \tau_{0}\lb({\xi_t(\bb)\atop \zx}\rb)$,
where $\tau_0$ is as defined above, and where we can treat
$\xi_t$ as a function of $\bb$ (because $\xi_t$ is 
$\SigAlg_*$-measurable).

This holds for any $\bb\in\bB^0$, so Claim 1.0 implies that 
$\mu^0_{|\SigAlg^*}[\bU_\zx^{t+1}]$ is the function
$\bD^{0}\ni
\ba\mapsto \tau_{0}\lb({\xi_t(\ba)\atop \zx}\rb)\in\CC{0,1}$.
But this function is  $\SigAlg\{t\}$-measurable (because 
$\xi_t$ is $\SigAlg\{t\}$-measurable by definition), so it
is also $\mu^0_{|\SigAlg\{t\}}[\bU_\zx^{t+1}]$.
\esubclaimprf

\subclaim{For any $\zx\in\sX$, \ 
$\mu^{(-\!1)}_{|\SigAlg^*}[\bU^{t+1}_\zx] = \mu^{(-\!1)}_{|\SigAlg\{t\}}[\bU^{t+1}_\zx] =
\tau_{-1}\lb({\xi_t\atop \zx}\rb)$, where

\centerline{\footnotesize $\D
\tau_{-1}\lb({{\hspace{1.1em}  \zl^t_2,\zl^t_1;\zd^t_0,\zd^t_1;\zr^t_1,\zr^t_2}
\atop {\hspace{-1.1em}\zl_2,\zl_1;\zd_0,\zd_1;\zr_1,\zr_2}}\rb)
\ \ := \ \
\choice{\D\frac{1}{(P_\bL)^2} & \If&
\mbox{\footnotesize
$\begin{array}[t]{l} 
\zl_2\in\sP_\bL(\zl_1), \ \ 
\zl_{1}\in\phi(\bL,\zl_{2}^t,\zl_1^t),  \\
\zd_0= \phi(\zl_{2}^t,\zl_1^t,\zd^t_0), \ \
\zd_1 = \phi(\zl_{1}^t,\zd^t_0,\zd^t_1), \\
\zr_1= \phi(\zd^t_0,\zd^t_1, \zr_{1}^t),
\And \zr_2 = \phi(\zd^t_1, \zr_{1}^t,\zr_2^t);
\end{array}$} \\
0 &&\mbox{otherwise.}}$}}
\bsubclaimprf
  Figure \ref{fig:diffuse2}(B) shows
how the values of $(\zd^{t+1}_0,\zd^{t+1}_1; \ \zr_{1}^{t+1},\zr_2^{t+1})$ are
determined by the data in $\SigAlg\{t\}$, because
$\zd^{t+1}_0 = \phi(\zl_{2}^t,\zl_1^t,\zd^t_0)$, \
$\zd^{t+1}_1 = \phi(\zl_{1}^t,\zd^t_0,\zd^t_1)$, \
$\zr_{1}^{t+1} = \phi(\zd^t_0,\zd^t_1, \zr_{1}^t)$,
and $\zr_{2}^{t+1} = \phi(\zd^t_1, \zr_{1}^t,\zr_2^t)$.
However, $\zl_{2}^{t+1} = \phi(\zl_{4}^t,\zl_{3}^t,\zl_{2}^t)$ 
and $\zl^{t+1}_{1} = \phi(\zl_{3}^t,\zl_{2}^t,\zl^t_{1})$
not determined, even  by $\SigAlg^*$.
Instead, for any fixed $z\in\Zahl$ and $\bb\in\bB_z^{(-\!1)}$,
let 
\[
\sL^0_\bb \ \ := \ \ 
\set{(\zl_{z\!-\!t\!-\!4},\zl_{z\!-\!t\!-\!3})\in\sA^2}{
\zl_{z\!-\!t\!-\!3}\in\sP_\bL(\zl^0_{z\!-\!t\!-\!2}) \And
 \zl_{z\!-\!t\!-\!4}\in\sP_\bL(\zl_{z\!-\!t\!-\!3})}; 
\]
Then $\#(\sL^0_\bb)=(P_\bL)^2$, because $\bL$ is left-regular. 
Let $\mu_{|\bb}$ be the conditional measure on $\bD^{(-\!1)}$
given $\bb$, and let $\lam^0_\bb:=
\pr{\{z\!-\!t\!-\!4,z\!-\!t\!-\!3\}}(\mu_{|\bb})\in\sM(\sL^0_\bb)$;
then $\lam^0_\bb$ is the uniform probability measure assigning
${1}/{(P_\bL)^2}$ to each element of $\sL^0_\bb$, because $\lam$ is the
Parry measure on $\bL$.

For any $(\zl_2,\zl_1)\in\bL_2$, let
\[
\sL(\zl_2,\zl_1)\quad:=\quad\set{(\zl'_2,\zl'_1)\in\sA^2}{
\zl'_1\in\phi(\bL,\zl_{2},\zl_1) \And \zl_2'\in\sP_\bL(\zl'_1)}.
\]
Then $\#\lb[\sL(\zl_2,\zl_1)\rb]=(P_\bL)^2$, 
because $\bL$ is left-regular and $\Phi$-left-resolving.
Let $\lam(\zl_2,\zl_1)\in\sM\lb[\sL(\zl_2,\zl_1)\rb]$
 be the uniform probability measure assigning
${1}/{(P_\bL)^2}$ to each element of $\sL(\zl_2,\zl_1)$.
Note that any $\bb\in\bB^{(-\!1)}$ completely determines the values of
$\zl_{2}^{t}$ and $\zl_{1}^{t}$ (because these
are $\SigAlg^*$-measurable functions).  
Let $\sL^{t+1}_\bb:=\sL\lb(\zl_{2}^{t}(\bb),\zl_{1}^{t}(\bb)\rb)$ and
$\lam^{t+1}_\bb:=\eta\lb(\zl_{2}^{t}(\bb),\zl_{1}^{t}(\bb)\rb)$.

\begin{sloppypar}
There is a function $\Phi_\bb:\sL^0_\bb\into\sL^{t+1}_\bb$ such that
$(\zl_{2}^{t+1},\zl_{1}^{t+1})
\ = \ \Phi_\bb(\zl_{z\!-\!t\!-\!4}^0,\zl_{z\!-\!t\!-\!3}^0)$, and
$\Phi_\bb$ is bijective because $\bL$ is
$\Phi$-left-resolving.
Thus, the $\mu^{(-\!1)}_{|\bb}$-conditional probability distribution of
$(\zl_{2}^{t+1},\zl_{1}^{t+1})$ is the measure $\Phi_\bb\lb[
\pr{\{z\!-\!t\!-\!4,z\!-\!t\!-\!3\}}(\mu_{|\bb})\rb]
=\Phi_\bb(\lam^0_\bb)\ \eeequals{(*)} \ \lam^{t+1}_\bb$. 
Here $(*)$ is because $\Phi_\bb$ is bijective,
while $\lam^0_\bb$ and $\lam^{t+1}_\bb$ are both uniform measures
on sets of $(P_\bL)^2$ elements.
Thus,
$\mu^{(-\!1)}_{|\bb}[\bU^{t+1}_\zx] =  \tau_{-\!1}\lb({\xi_t(\bb)\atop \zx}\rb)$,
where $\tau_{-\!1}$ is as defined above, and where we can again treat
$\xi_t$ as a function of $\bb$.
\end{sloppypar}

This holds for any $\bb\in\bB^{(-1)}$, so Claim 1.0 implies that 
$\mu^{(-\!1)}_{|\SigAlg^*}[\bU_\zx^{t+1}]$ is the function
$\bD^{(-\!1)}\ni
\ba\mapsto \tau_{-\!1}\lb({\xi_t(\ba)\atop \zx}\rb)\in\CC{0,1}$.
But this function is  $\SigAlg\{t\}$-measurable (because 
$\xi_t$ is $\SigAlg\{t\}$-measurable), so it
is also $\mu^{(-\!1)}_{|\SigAlg\{t\}}[\bU_\zx^{t+1}]$.
\ignore{
 This holds for any $\bb\in\bB^{(-\!1)}$, so Claim 1.0 implies that the
$\mu^{(-\!1)}_{|\SigAlg^*}$-conditional probability distribution of
$(\zl_{2}^{t+1},\zl_{1}^{t+1})$ is the function
$\D
\sum_{\bb\in\bB^{(-\!1)}}
\mu^{(-\!1)}[\bb]\cdot \hlam^{t+1}_\bb$,
where $\hlam^{t+1}_\bb :=\chr{\sL^{t+1}_\bb}$.
  But this function is $\SigAlg\{t\}$-measurable (because it is
a sum of $\SigAlg\{t\}$-measurable characteristic functions),
so it is also the $\mu^{(-\!1)}_{|\SigAlg\{t\}}$-conditional probability distribution of
$(\zl_{2}^{t+1},\zl_{1}^{t+1})$.  The claim follows.}
\esubclaimprf

\subclaim{For any $\zx\in\sX$, \
$\mu^1_{|\SigAlg^*}[\bU^{t+1}_\zx]  =  \mu^1_{|\SigAlg\{t\}}[\bU^{t+1}_\zx]
 =  \tau_{1}\lb({\xi_t\atop \zx}\rb)$, where 
$\tau_1:\sX\into\sM(\sX)$ is defined similarly to $\tau_{-1}$.}
\bsubclaimprf  Figure \ref{fig:diffuse2}(C) shows
how  the values of $(\zl^{t+1}_2,\zl^{t+1}_1,\zd^{t+1}_0,\zd^{t+1}_1)$ are
determined by the data in $\SigAlg\{t\}$, because
$\zl_{2}^{t+1} = \phi(\zl^t_2,\zl^t_1,\zd^t_0)$,
$\zl_{1}^{t+1} = \phi(\zl^t_1,\zd^t_0,\zd^t_1)$,
$\zd^{t+1}_0 = \phi(\zd^t_0,\zd^t_1,\zr^t_1)$,
and
$\zd^{t+1}_1 = \phi(\zd^t_1,\zr^t_1,\zr^t_2)$.
However, $\zr_{1}^{t+1} = \phi(\zr_{1}^t,\zr^t_{2},\zr^t_3)$
and $\zr^{t+1}_{2} = \phi(\zr_{2}^t,\zr^t_{3},\zr^t_4)$
not determined, even by $\SigAlg^*$.
Now proceed as in Claim 1.2, but
replace $\zl_{2}^{t+1}$ by $\zr_{2}^{t+1}$, \
$\zl_{1}^{t+1}$ by $\zr_{1}^{t+1}$, \
$\bL$ with $\bR$, \  $\lam$ with $\rho$,  and `left-resolving' with 
`right-resolving'.
\esubclaimprf

Finally, define $\tau:\sX\into\sM(\sX)$ by $\tau(\zy,\zx) := 
\tau_{\vV(\zy)}\lb({\zy\atop \zx}\rb)$,
where $\vV(\zy)\in\{-1,0,1\}$, and where $\tau_{0}$ and $\tau_{\pm1}$ are
defined as in  Claims 1.1 to 1.3.  Then
$\mu_{|\SigAlg^*}[\bU^{t+1}_\zx]  =  \mu_{|\SigAlg\{t\}}[\bU^{t+1}_\zx]
 =  \tau(\xi_t; \zx)$.
\eclaimprf

\Claim{For any $\zx\in\sX$ and $t\in\Natur$, \
$\mu_{|\SigAlg\CC{0...t}}[\bU^{t+1}_\zx] \ = \ \mu_{|\SigAlg^*}[\bU^{t+1}_\zx]$. 
}
\bclaimprf
Claim 1 implies that
$\mu_{|\SigAlg^*}[\bU^{t+1}_\zx]$ is actually 
a $\SigAlg\{t\}$-measurable function,
because it is equal to $\mu_{|\SigAlg\{t\}}[\bU^{t+1}_\zx]$.  But this means
that $\mu_{|\SigAlg^*}[\bU^{t+1}_\zx]$ is also $\SigAlg\CC{0...t}$-measurable,
because $\SigAlg\{t\} \subseteq \SigAlg\CC{0...t}$.  Also, for any $\bC\in \SigAlg\CC{0...t}$,
we have $\int_\bC \mu_{|\SigAlg^*}[\bU^{t+1}_\zx] \ d\mu = 
\mu[\bC\intsct\bU^{t+1}_\zx]$, because $\bC\in\SigAlg^*$, because 
$\SigAlg\CC{0...t}\subseteq\SigAlg^*$.
  But $\mu_{|\SigAlg\CC{0...t}}[\bU^{t+1}_\zx]$
is the unique $\SigAlg\CC{0...t}$-measurable function with this property
(by definition); hence $\mu_{|\SigAlg^*}[\bU^{t+1}_\zx] =  \mu_{|\SigAlg\CC{0...t}}[\bU^{t+1}_\zx]$.
\eclaimprf
For any  $\zx\in\sX$ and $t\in\Natur$,
we conclude
$\mu_{|\SigAlg\CC{0...t}}[\bU^{t+1}_\zx] $ $
\ \ \eeequals{(*)} \ \ $ $
(\mu_{|\SigAlg^*})[\bU^{t+1}_\zx] $ $
\ \ \eeequals{(\dagger)}\ \ $ $
\mu_{|\SigAlg\{t\}}[\bU^{t+1}_\zx] $ $ \ \ \eeequals{(\dagger)}\ \ $ $ \tau(\xi_t,\zx)$,
where $(*)$ is Claim 2 and $(\dagger)$ is Claim 1.
Thus, eqn.(\ref{markov.eqn}) is satisfied. 
\ethmprf

\NRemark{\label{rem:defect.random.walk}
(a) Suppose the conditions of Theorem \ref{defect.random.walk} are
satisfied.  Then the measure $\del$ can always be chosen so that the
Markov measure $\barmu$ is shift-invariant (because every finite-state
Markov chain has a stationary measure).  The $\shift{}$-ergodic
components of $\barmu$ are then the stochastic analogs of the
{\dfn particle types}  of  Definition \ref{def:type}.

The {\dfn drift velocity} of $\omg$ is the expected value 
$\Vdrift(\omg):=\D\sum_{\fv\in\dV} \nu[\fv]\cdot \fv$. 
If $\barmu$ is $\shift{}$-ergodic 
(i.e. $\barmu$ corresponds to a single particle type),  then
for $\forall_\omg \ \bz\in\Zahl^\Natur$,
\ $\D\lim_{n\goto\oo} (z_n/n)  =  \Vdrift(\omg)$
(by the Birkhoff Ergodic Theorem).  Thus,
$\Vdrift(\omg)$ is the long-term average velocity of particles of type $\barmu$.
For instance, the Markov chain in Example \ref{X:diffuse} has
one ergodic component (i.e. one particle type), with $\Vdrift=0$.

\begin{sloppypar}
(b) A special case of Theorem \ref{defect.random.walk} was previously
proved in \cite[Thm.2.1.1]{Elo93b},
for when $\bL=\LZ$ and $\bR=\RZ$ are permutative
subalphabets for $\Phi$ [see Example \ref{X:resolve}(a)] and $W=0$.
In this case we must have
$\sR\neq\sL$ for an $(\LZ,\RZ)$-defect to be meaningful.
We recommend \cite{Elo93a,Elo93b} for further interesting examples 
of diffusive defect dynamics, as well as analysis of their
drift and variance.
These methods were extended to defect ensembles
in \cite{Elo94}, and to the pseudorandom motion of domain
boundaries in two-dimensional boolean CA \cite{Elo95}. 
\end{sloppypar}

(c) Empirically, the large $\alp$ defect particle of ECA\#54 
[see Figure \ref{fig:54.defects}$(\alp)$] also performs a random walk,
as can perhaps be seen  in Figure \ref{fig:defect.intro}(A).
However, this motion is not due to the mechanism of
Theorem \ref{defect.random.walk}, because $\alp$ belongs to the
`ballistic' regime of \S\ref{S:ballistic}, not the `diffusive' regime.
Instead, the meandering is due to interactions
with neighbouring $\alp$ particles, mediated by a complex exchange
of the tiny $\gam^\pm$ particles of Figure \ref{fig:54.defects}$(\gam^\pm)$.
See \cite[Fig.13(b)]{CrHa97}.
}

\Corollary{\label{defect.random.walk2}}
{
  Let $\Phi:\AZ\into\AZ$ be a CA and fix $p,q\in\Natur$.  Suppose that
either
\bitem
  \item[]{\bf[i]} \
$\bL\subseteq\Fix{\Phi^p,\shift{q}}$ and $\lam\in\sM(\bL)$, while 
$\bR\subseteq\AZ$ is a right-resolving, right-regular Markov subshift
with Parry measure $\rho\in\sM(\bR)$.

  \item[or] {\bf[ii]} \  $\bL\subseteq\AZ$ is a left-resolving,
left-regular Markov subshift with Parry measure
 $\lam\in\sM(\bL)$, while
$\bR\subseteq\Fix{\Phi^p,\shift{q}}$ and $\rho\in\sM(\bR)$.
\eitem
Let $W\in\Natur$, let $\del$ be any probability measure on
$\sD:=\sA^W$, and let $\mu:=\lam\tensor\del\tensor\rho
\in\sM(\bD^{W,0}_{\bL,\bR})$.
Define $\zeta:\bD^W_{\bL,\bR}\into\Zahl^\Natur$
by $\zeta(\ba):=(z_0,z_p,z_{2p},\ldots)$, where, for all $t\in\Natur$,
\ $z_{tp}\in\Zahl$ is as in 
eqn.{\rm(\ref{variable.width.defect.2})}.
Then $\omg:=\zeta(\mu)\in\sM(\Zahl^\Natur)$ is a random walk.
}
\bthmprf ({\em Case} {\bf[i]})
By using the $q$th higher power representation of $\AZ$
[see Remark \ref{defect.particle.remarks}(c)]
and replacing $\Phi$ with $\Phi^p$, we can assume that $q=p=1$;
i.e. that $\bL\subseteq \Fix{\Phi,\shift{}}$.  Thus
$\bL=\{\zl^\oo\}_{\zl\in\sL}$, where $\sL\subseteq\sA$ is some
subalphabet, and, for each $\zl\in\sL$, the point $\zl^\oo:=
[\ldots \zl\zl\zl\ldots]$ is $\Phi$-fixed.
Thus $\lam= \D\sum_{\zl\in\sL} c_\zl \chr{\zl}$,
where, for each $\zl\in\sL$, $\chr{\zl}$ is the point-mass
on $\zl^\oo$, and $c_\zl\in\CC{0,1}$ is a constant.
Thus, $\mu =  \D\sum_{\zl\in\sL} c_\zl \mu_\zl$, where
$\mu_\zl:=\chr{\zl}\tensor\del\tensor\rho$.  It suffices to
prove that each of the measures $\mu_\zl$ induces a random walk.
Hence, assume that $\lam=\chr{\zl}$ for some $\zl\in\sL$,
and redefine $\sL:=\{\zl\}$ and $\bL:=\{\zl^\oo\}$.
Then $\bL$ is a left-regular and left-resolving
subshift [see Example \ref{X:resolve}(e)]. 
Now apply Theorem \ref{defect.random.walk}.
The proof of {\em Case} {\bf[ii]} is analogous.
\ethmprf

\section{The Turing Regime and Pushdown Regimes \label{S:turing}}

  Recall that a Turing machine \cite[\S7.2]{HopcroftUllman}
consists of a `head' which
deterministically moves back and forth along a `tape', reading and
writing symbols from some alphabet.  To be precise, let $\sT$ 
be a finite set.  
A (classical) {\dfn Turing Machine} with {\dfn tape alphabet}
$\sT$ is a quadruple $(\sD,\tau,\Ups,\vV)$, where
$\sD$ is a finite set (called the {\dfn head state domain}), 
$\tau:\sT\x\sD\into\sT$ is a {\dfn tape rule},
 $\Ups:\sT\x\sD\into\sD$ is an {\dfn update rule},
 and $\vV:\sT\x\sD\into\{-1,0,1\}$ is a {\dfn velocity rule}.
The {\dfn machine statespace} of the Turing machine is $\sT^\Zahl\x\sD\x\Zahl$.
If the machine is in state $(\bt,\zd,z)\in\sT^\Zahl\x\sD\x\Zahl$, this
means that the tape  currently has symbol string $\bt$, the head is
at position $z$ on the tape, and the head has state description $\zd$. 
If $\bt:=[\ldots \zt_{z-1} \ \zt_z \ \zt_{z+1}\ldots]$, then
define $\bt' :=  [\ldots \zt_{z-1} \ \tau(\zt_z,\zd) \ \zt_{z+1}\ldots]$.
The dynamics of the machine is the map $\Theta:\sT^\Zahl\x\sD\x\Zahl\into
\sT^\Zahl\x\sD\x\Zahl$ defined:
\[
  \Theta(\bt,\zd,z) \quad:=\quad \lb(\maketall \bt', \ \Ups(\zt_z,\zd), \ z+\vV(\zt_z,\zd)\rb).
\]
 We will generalize this definition in two ways.  First,
we will imagine that the head lies {\em between} two tape symbols,
rather than {\em over} a tape symbol.  The head can read the
two symbols to its left and two symbols to its right, and
can overwrite the symbol immediately left or right.
Second, we require that there are 
Markov subshifts $\bL,\bR\subseteq\AZ$ such that
the symbol sequence on the left half of the
tape lies $\bL^-$, while the right half
lies in $\bR^+$.  The machine must write new symbols so as to respect the constraints of these subshifts.  

\newcommand{\tauL}{\tau_{\scriptscriptstyle L}}
\newcommand{\tauC}{\tau_{\scriptscriptstyle C}}
\newcommand{\tauR}{\tau_{\scriptscriptstyle R}}

\newcommand{\betL}{\bet_{\scriptscriptstyle L}}
\newcommand{\betC}{\bet_{\scriptscriptstyle C}}
\newcommand{\betR}{\bet_{\scriptscriptstyle R}}

Formally, an {\dfn $(\bL,\bR)$-Turing machine} is
a sextuple $(\sD,\tauL,\tauC,\tauR,\Ups,\vV)$, where
$\sD$ is a finite set, 
$\tauL:\sA^2\x\sD\into\sA$, \
$\tauC:\sA\x\sD\x\sA\into\sA$, \
$\tauR:\sD\x\sA^2\into\sA$, \
 \ $\Ups:\sA^2\x\sD\x\sA^2\into\sD$,
 and $\vV:\sA\x\sD\x\sA\into\{-1,0,1\}$.
The {statespace} of the Turing machine is $\bL^-\x\sD\x\bR^+\x\Zahl$.
If the machine is in state $(\bl,\zd,\br,z)\in\bL^-\x\sD\x\bR^+\x\Zahl$, this
means that the tape  currently has symbol string $(\bl \ | \ \br)$, 
where the head (indicated by `$|$') is at position
  $z+\frac{1}{2}$ on the tape, and the head has state description $\zd$. 
The machine dynamical system  $\Theta:\bL^-\x\sD\x\bR^+\x\Zahl\into
\bL^-\x\sD\x\bR^+\x\Zahl$ is defined
by $\Theta(\bl,\zd,\br,z) := \   ( \bl', \ \zd', \ \br', z')$,
where $\zd' \ : = \ \Ups(\zl_2,\zl_1,\zd,\zr_1,\zr_2)$ and
$z' \ := \   z+\vV(\zl_1,\zd,\zr_1)$ and
\beqn
\label{LR.turing.defn}
\hspace{-2em}(\bl' \ | \ \br')  :=  \lb\{
\begin{array}{rclcl}
(\ldots,\zl_{4},\zl_{3},\zl_{2} \ \ \  | &  \zr_0',\! &\zr_1',\zr_2,\zr_3,\zr_4,\ldots)  &\If & \vV(\zl_1,\zd,\zr_1)=-1;\\
(\ldots,\zl_{4},\zl_{3},\zl_{2},\zl_1' & |  &  \zr_1',\zr_2,\zr_3,\zr_4,\ldots)
& \If &  \vV(\zl_1,\zd,\zr_1)=0;\\
(\ldots,\zl_{4},\zl_{3},\zl_{2},\zl_1',\! & \zl_0' & |  \ \ \ \zr_2,\zr_3,\zr_4,\ldots)  & \If & \vV(\zl_1,\zd,\zr_1)=+1. \end{array}\rb.\\
\eeqn
\[
\begin{array}{rclrrcll}
\mbox{Here,} \ 
\zl_1' &:=& \tauL(\zl_2,\zl_1;\zd) & \mbox{is such that} &
 (\zl_2,\zl_1')&\in&\bL_2; \\
\zr_1' &:=& \tauC(\zd;\zr_1,\zr_2) & \mbox{is such that} & (\zr_1',\zr_2)&\in&\bR_2; \\
\And \zl_0' &:=& \tauC(\zl_1;\zd;\zr_1) & \mbox{is such that} &
(\zl_1',\zl_0')&\in&\bL_2, & \If  \vV(\zl_1,\zd,\zr_1)=+1, \\ 
\mbox{whereas} \ \zr_0' &:=& \tauC(\zl_1;\zd;\zr_1) & \mbox{is such that} &
(\zr_0',\zr_1')&\in&\bR_2 & \If  \vV(\zl_1,\zd,\zr_1)=-1. 
\end{array}
\]
(If  $\vV(\zl_1,\zd,\zr_1)=0$, then the value of $ \tauC(\zl_1;\zd;\zr_1)$
is discarded, so it is irrelevant). Finally,
\bitem
  \item $\Ups(\zl_2,\zl_1,\zd,\zr_1,\zr_2)$ depends only on $(\zl_2,\zl_1,\zd)$
if $\vV(\zl_1,\zd,\zr_1)=-1$.
  \item $\Ups(\zl_2,\zl_1,\zd,\zr_1,\zr_2)$ depends only on $(\zl_1,\zd,\zr_1)$
if $\vV(\zl_1,\zd,\zr_1)=0$.
  \item $\Ups(\zl_2,\zl_1,\zd,\zr_1,\zr_2)$ depends only on $(\zd,\zr_1,\zr_2)$
if $\vV(\zl_1,\zd,\zr_1)=+1$.
\eitem

\Proposition{\label{Turing.CA}}
{
Let $\bL,\bR\subset\AZ$ be Markov subshifts.   Let $W\in\Natur$
and let $\sD:=\sA^W$.  Let $\hbL,\hbR\subset\sD^\Zahl$ be the
$W$th higher power representations of $\bL$ and $\bR$.
\bthmlist
  \item Let $\Phi:\AZ\into\AZ$ be a CA with
$\bL,\bR\subseteq\Fix{\Phi}$.  
Then the dynamical system $(\bD^W_{\bL,\bR},\Phi)$ is isomorphic
to an $(\hbL,\hbR)$-Turing machine $(\sD,\tau,\Ups,\vV)$.

  \item Conversely, given any $(\hbL,\hbR)$-Turing machine $(\sD,\tau,\Ups,\vV)$, there is a CA $\Phi:\AZ\into\AZ$,
with $\bL,\bR\subseteq\Fix{\Phi}$, 
such that  $(\bD^W_{\bL,\bR},\Phi)$ is isomorphic to 
$(\sD,\tau,\Ups,\vV)$.
\ethmlist
}
{\em Proof idea:}  The defect acts like the Turing machine head.
The application of $\Phi$ changes the head state, and can
can also modify the adjacent symbols on the $(\bL,\bR)$-tape.  However,
just as in a Turing machine, the
more distant tape symbols remain unchanged, because $\bL$ and $\bR$
are $\Phi$-fixed.

\bthmprf {\bf(a)} By passing to the $W$th higher power recoding, 
and replacing $\bL$ with $\hbL$ and $\bR$ with $\hbR$
and $\sA$ with $\hsA:=\sA^W$, Remark \ref{defect.particle.remarks}(c)
allows us to assume that $\sD=\hsA^2$ and that
$\Upsilon:\hsA\x\sD\x\hsA \into \sD$ and 
$\vV:\hsA\x\sD\x\hsA \into \{-1,0,1\}$
in eqn.(\ref{defect.finite.state.update}).  To simplify notation,
we will suppress the `hats'. 
 Define $\Psi:\bD^2_{\bL,\bR}\into\bL^-\x\sD\x\bR^+\x\Zahl$
so that, if $\ba^t$ is as in eqn.(\ref{variable.width.defect.2}), then
$\Psi(\ba^t):=(\bl,\bd,\br;z)$, where
\[
\begin{array}{rclcrcl}
\bl&:=&[\ldots,\zl^t_3,\zl^t_2,\zl^t_1]\in\bL^-, 
&& 
\bd&:=&[d^t_{0}, d^t_1]\in\sD=\sA^2,\\
\br&:=&[\zr^t_1,\zr^t_2,\zr^t_3,\ldots]\in\bR^+,
&\And& z&:=&z_t \ \in \ \Zahl.
\end{array}
\]
Let $\Ups$ and $\vV$ be as in
 eqn.(\ref{defect.finite.state.update}).  If $\ba^{t+1}:=\Phi(\ba^t)$, then
$\Psi(\ba^{t+1}) = (\bl',\bd',\br';z')$, where $\bd':=
[d^{t+1}_{0}, d^{t+1}_{1}] = \Ups(\zl_1,\bd,\zr_1)$, 
and where $\bl'$ and $\br'$ are as in eqn.(\ref{LR.turing.defn}),
with
\beq
\ \tauL(\zl_2,\zl_1,\bd) &:=& \Phi(\zl_2,\zl_1,d^t_{0});\\
\tauC(\zl_1,\bd,\zr_1) &:=&
\choice{ \Phi(\zl_1,d^t_{0},d^t_{1}) & \If & \vV(\zl_1,\bd,\zr_1)=+1; \\
\Phi(d^t_{0},d^t_{1},\zr_1) & \If & \vV(\zl_1,\bd,\zr_1)=-1;\\
 \mbox{irrelevant} & \If & \vV(\zl_1,\bd,\zr_1)=0.} \\
\And \tauR(\bd,\zr_1,\zr_2) &:=& \Phi(d^t_{1},\zr_1,\zr_2). 
\eeq
{\bf(b)} is a straightforward generalization of the method of Lindgren and Nordahl
\cite{LindgrenNordahl} for simulating a classical Turing machine with
a cellular automaton.
\ethmprf

  Proposition \ref{Turing.CA} applies even when $\bL$ and $\bR$ are
$\shift{}$-periodic subshifts, but in this case it isn't very interesting,
because an $(\bL,\bR)$-admissible `tape'
can't encode any information, so the resulting Turing machine is
rather trivial, and is described in \S\ref{S:ballistic}.
To perform useful computation, we need
$\bL$ and $\bR$ to have nonzero entropy.
 If $\gB\subset\AZ$ is a subshift, then the {\dfn topological entropy}
of $\gB$ is defined
\[
  h(\gB,\shift{}) \quad:=\quad 
\lim_{N\goto\oo} \frac{\log_2(\#\gB_{\CO{0...N}})}{n}.
\]
If $\gB$ is a subshift of finite type, then $h(\gB,\shift{})>0$ iff
$\gB$ is not $\shift{}$-periodic. In particular, if
$\gB$ is a Markov subshift defined by a digraph on the vertex set
$\sA$, then $h(\gB,\shift{})>0$ iff this digraph is not just a disjoint
union of cycles.  Equivalently, there is
a {\dfn choice point} vertex $c\in\sA$,  
meaning that $c$ belongs to at least two distinct
cycles.  See \cite[Ch.4]{LindMarcus}, \cite[\S1.4]{Kitchens} or
\cite[\S3.6.2]{KurkaBook}.  A $\gB$-admissible sequence $\bb$ can then
encode nontrivial information, because for every $z\in\Zahl$ with
$b_z=c$, there are at least two $\gB$-admissible possibilities for
$b_{z+1}$, and a choice between these encodes at least one bit of
information.

A pushdown automaton \cite[\S5.2]{HopcroftUllman}
is a finite automaton augmented with a 
`stack' or `last in, first out' (`LIFO') memory model.
To be precise, a {\dfn pushdown automaton} is a septuple 
$(\sI,\sD,\sO,\sT;\Ups,\Omg,\Sigma)$, where
$\sI$, $\sD$, and $\sO$ are a finite input space, state domain,
and output space, respectively (as in a finite automaton), and
$\sT$ is a finite {\dfn stack alphabet}.  Now
$\Ups:\sI\x\sT\x\sD\into\sD$ is the {\dfn update rule},
$\Omg:\sI\x\sT\x\sD\into\sO$ is the {\dfn output rule},
and $\Sigma:\sI\x\sT\x\sD\into\sT\disj\{\emptyset,\triangle\}$ is a {\dfn stack rule}.
The {\dfn machine statespace} of the pushdown automaton is 
$\sD\x\sT^\Natur$. 
The machine behaviour is defined by the map $\Theta:
\sI\x\sD\x\sT^\Natur\into\sD\x\sT^\Natur\x\sO$ defined
$\Theta(\zi,\zd,\bt) \ :=  (\zd',\bt',\zo)$,
where $\zd'=\Ups(\zi,\zd,\zt_0)$, $\zo:=\Omg(\zi,\zd,\zt_0)$, and where
\[
\bt' \quad := \quad 
\choice{ (\zt_1,\zt_2,\zt_3,\ldots) & \If & \Sigma(\zi,\zd,\zt_0)=\triangle
\quad \mbox{\footnotesize {(i.e. `pop' the symbol $\zt_0$ off the stack);}} \\
         (\zt_0,\zt_1,\zt_2,\ldots) & \If & \Sigma(\zi,\zd,\zt_0)=\emptyset
\quad \mbox{\footnotesize {(i.e. do not touch the stack);}} \\
         (\zt',\zt_0,\zt_1,\ldots) & \If & \Sigma(\zi,\zd,\zt_0)=\zt'
\quad \mbox{\footnotesize {(i.e. `push' the symbol $\zt'$ onto the stack).}}} 
\]
An {\dfn autonomous finite automaton} is a finite
automaton with no input or output; i.e. a dynamical system
$\Ups:\sD\into\sD$ where $\sD$ is a finite set.
Similarly, an {\dfn autonomous pushdown automaton} (APDA) is a pushdown
automaton with no input or output; i.e. $\sI=\emptyset=\sO$.  Thus, the
function $\Omg$ is trivial, $\Ups:\sT\x\sD\into\sD$, and
$\Sigma:\sT\x\sD\into\sT\disj\{\emptyset,\triangle\}$.  The APDA's
future behaviour is entirely determined by the initial stack state.
It is easy to see that an APDA is equivalent to an $(\bL,\bR)$-Turing
machine where $\bR=\AZ$ and $\bL=\{\bar0\}$ where $0$ is some `null'
symbol.   Thus, we will treat it as such.

Let $\fM$ and $\fM'$ be two machine-classes.  We write $\fM\preceq\fM'$
if any machine in $\fM$ can
be simulated by one in $\fM'$  (possibly not in real time).
We say that $\fM$ and $\fM'$  are {\dfn computationally equivalent} 
(and write $\fM\approx\fM'$)
if $\fM\preceq\fM'$ and $\fM'\preceq\fM$.
Let $\Turing$ be the class of (classical) Turing machines,
and let $\Turing_{\bL,\bR}$ be the class of $(\bL,\bR)$-Turing machines.
Let $\APDA$ be the class of autonomous pushdown automata, and let
$\AFA$ be the class of autonomous finite automata.

\Proposition{\label{machine.equivalence}}
{
Let $\bL,\bR\subset\AZ$ be Markov subshifts.
\bthmlist
\item If
$h(\bL,\shift{})>0$  and $h(\bR,\shift{})>0$,
then $\Turing_{\bL,\bR} \approx \Turing$.

\item If
$h(\bL,\shift{})>0=h(\bR,\shift{})$,
or 
$h(\bL,\shift{})=0<h(\bR,\shift{})$,
then $\Turing_{\bL,\bR} \approx \APDA$.

\item If
$h(\bL,\shift{})=0=h(\bR,\shift{})$,
then $\Turing_{\bL,\bR} \approx \AFA$.
\ethmlist
}
\bthmprf
 By a {\dfn cycle of length $P$} in $\bR$, we mean a word
$\bc=(c_1,\ldots,c_{P})\in\bR_{P}$, such that $(c_P,c_1)\in\bR_2$;
hence the infinite sequence $[\ldots \bc\bc\bc\ldots]$ is $\bR$-admissible.

\Claim{Suppose $h(\bR,\shift{})>0$.
Then there is some $P\in\Natur$ and some $c\in\sA$ such that
$c$ begins two different cycles $\bc_0$ and $\bc_1$ 
in $\bR$, both of length $P$.}
\bclaimprf
$h(\bR,\shift{})>0$, so
there is some $c\in\sA$
which belongs to two different cycles in $\bR$; say 
$\bb_0=(b^0_1,b^0_2,\ldots,b^0_{Q_0})$
and $\bb^1=(b^1_1,b^1_2,\ldots,b^1_{Q_1})$, where $Q_0,Q_1\in\Natur$
and $b^0_1=c=b^1_1$.  Let $P:=\lcm(Q_0,Q_1)$.
Let $\bc_0$ (resp. $\bc_1$) 
 be the cycle obtained by chaining together $P/Q_0$ copies of $\bb_0$
(resp. $P/Q_1$ copies of $\bb_1$).
Then $\bc_0$ and $\bc_1$ are distinct cycles of length $P$, both starting with $c$.
\eclaimprf

\Claim{Suppose $h(\bR,\shift{})>0$ and $h(\bL,\shift{})>0$.
Then there is some $P\in\Natur$ and some $\zr,\zl\in\sA$ such that
$\zr$ begins two different cycles $\br_0$ and $\br_1$ 
in $\bR$, and
$\zl$ begins two different cycles $\bl_0$ and $\bl_1$ 
in $\bL$,  with all four cycles having length $P$.}
\bclaimprf
  Claim 1 yields two cycles $\bc^R_0,\bc^R_1$ in $\bR$, say of length $P_R$,
beginning with the same symbol, say $\zr$.
The same argument also yields two cycles $\bc^L_0,\bc^L_1$
in $\bL$, say of length $P_L$, 
beginning with the same symbol, say $\zl$.  Let
$P:=\lcm(P_R,P_L)$.
Let $\br_0$ (resp. $\br_1$) 
 be obtained by chaining together $P/P_R$ copies of $\bc^R_0$
(resp. $\bc^R_1$).
Let $\bl_0$ (resp. $\bl_1$) 
 be obtained by chaining together $P/P_L$ copies of $\bc^L_0$
(resp. $\bc^L_1$).
\eclaimprf

Let $\bR_1^+\subseteq\bR^+$ be the set of all right-infinite sequences
made by concatenating copies of $\br_0$ and $\br_1$, and let
$\sT:=\{0,1\}$.  Define $\betR:\sT\into\{\br_0,\br_1\}$
by $\betR(\zt):=\br_\zt$ for $\zt=0,1$.
Define bijection
$\betR^\Natur:\sT^\Natur\into\bR_1^+$  by $\betR(\zt_1,\zt_2,\zt_3,\ldots):=
[\betR(\zt_1)\ \betR(\zt_2) \ \betR(\zt_3) \ \ldots]$. 
 Clearly, $\betR\circ\shift{} =
\shift{P}\circ\betR$.  In this way, we can encode any binary sequence
with an element of $\bR_1^+$.  Likewise, let $\bL_1^-\subseteq\bL^-$
be the set of all left-infinite sequences made from $\bl_0$ and
$\bl_1$, define $\betL:\sT\into\{\bl_0,\bl_1\}$ by
$\betL(\zt):=\bl_\zt$ for $\zt=0,1$,
and define bijection $\betL^\Natur:\sT^{-\Natur}\into\bL_1^-$ by
$\betR(\ldots,\zt_{-3},\zt_{-2},\zt_{-1}):=
[\ldots \ \betL(\zt_{-3}) \ \betL(\zt_{-2}) \ \betL(\zt_{-1})]$.  Clearly,
$\betL\circ\shift{-1} = \shift{-P}\circ\betL$.

{\bf(a)} \  ``$\Turing_{\bL,\bR} \preceq \Turing$'':
Every $(\bL,\bR)$-Turing machine is clearly also an
 $(\AZ,\AZ)$-Turing machine, which can clearly be simulated by a 
classical Turing machine with tape alphabet $\sA$.

``$\Turing \preceq \Turing_{\bL,\bR} $'': \ If
$\TM=(\sD,\tau,\Ups,\vV)$ is a classical Turing machine with tape
alphabet $\sT$, then one tape symbol $\zt_0$ lies directly `underneath'
the head.  In contrast, in an $(\bL,\bR)$-Turing machine
$\TM_*=(\sD_*,\tauL,\tauC,\tauR,\Ups_*,\vV_*)$, we only have tape
symbols to the left and right sides.  However, if the head state
domain $\sD_*$ is large enough, then the head of $\TM_*$ can temporarily
`remember' the value of the symbol $\zt_0$, even though $\zt_0$ is not
written anywhere on the tape.  So, let $\sD_*:=\sD\x\sT\x\bL_P\x\bR_P$.
If $\bt:=(\ldots,\zt_{-2},\zt_{-1},\zt_0,\zt_1,\zt_2,\ldots)$
and $\TM$ has headstate $\zd\in\sD$, then 
let $\bl:=\betL^\Natur(\ldots,\zt_{-3},\zt_{-2},\zt_{-1})$,
let $\br:=\betR^\Natur(\zt_1,\zt_2,\zt_3,\ldots)$,
and let $\TM_*$ have head state description $\zd_*:=(\zd,\zt_0,\bl_*,\br_*)\in\sD_*$.  The third and fourth entries of $\zd_*$ are input buffers; their values
(represented by $\bl_*$ and $\br_*$) are currently irrelevant.

Now, suppose $\TM$ moves {\em right} by one step, overwriting $\zt_0$ with
$\zt_0'$ and changing its head state to $\zd'$.
Then $\TM_*$ moves right by $P$ steps; 
during which time it reads $[\zr_{1},\ldots,\zr_{P}]$ and stores this
$P$-tuple in the fourth entry of $\zd_*$ (labelled `$\br_*$' above), while
writing the $P$ symbols of $\betL(\zt_0')$ to the tape.
Finally $\TM_*$ computes $\zt_1:=\betR^{-1}[\zr_{1},\ldots,\zr_{P}]$
and changes its headstate to $\zd'_*:=(\zd',\zt_1,\bl_*,\br'_*)$.
(where $\bl_*$ and $\br'_*$ are again irrelevant).

Suppose $\TM$ moves {\em left} by one step, overwriting $\zt_0$ with
$\zt_0'$ and changing its head state to $\zd'$.  Then $\TM_*$ moves
left by $P$ steps, during which time it reads
$[\zl_{P},\ldots,\zl_{1}]$ (in reverse order) and stores this $P$-tuple
in the third entry of $\zd_*$ (labelled `$\bl_*$' above), while writing the
$P$ symbols of $\betR(\zt_0')$ to the tape.  Finally $\TM_*$ computes
$\zt_{-1}:=\betL^{-1}[\zl_{P},\ldots,\zl_{1}]$ and changes its
headstate to $\zd'_*:=(\zd',\zt_{-1},\bl'_*,\br_*)$
(where $\bl'_*$ and $\br_*$ are again irrelevant).

Thus, the update rule $\Ups_*$ not only must emulate $\Ups$, but also must
implicitly compute $\betL^{-1}$
and $\betR^{-1}$.  Also, the tape rules $\tauL$ and $\tauC$ 
not only must emulate $\tau$, but also must
implicitly compute $\betL$;  likewise, the tape rules $\tauC$ and $\tauR$ must
implicitly compute $\betR$.

{\bf(b)} \ Suppose $h(\bL,\shift{})=0<h(\bR,\shift{})$
(the case ``$h(\bL,\shift{})>0=h(\bR,\shift{})$'' is analogous).

``$\APDA \preceq \Turing_{\bL,\bR} $'':
Let $\bc_0,\bc_1$ be as in  Claim 1, and
let $\bR_1^+\subseteq\bR^+$ be the set of all right-infinite sequences
made by concatenating copies of $\bc_0$ and $\bc_1$.
Now define a bijection $\betR^\Natur:\sT^\Natur\into \bR_1^+$,
and use $\betR^\Natur$ to build an $(\bL,\bR)$-Turing machine which
can emulate a given APDA with stack alphabet $\sT$, as
in part (a).  

 ``$\Turing_{\bL,\bR} \preceq \APDA$'':
 $\bL$ is $\shift{}$-periodic, so 
by passing to a higher power presentation, we can assume $\bL$ contains
only constant sequences.  At this point, any $(\bL,\bR)$-Turing machine
is clearly an APDA.

{\bf(c)} ``$\Turing_{\bL,\bR}\preceq \AFA$'': 
 If $h(\bL,\shift{})=0=h(\bR,\shift{})$, then both 
$\bL$ and $\bR$ are periodic,
so by passing to a higher power presentation, we can assume that
both $\bL$ and $\bR$
contain only constant sequences.  Thus,  the only computation performed by an
$(\bL,\bR)$-Turing machine is computation of the update rule
$\Ups:\sD\into\sD$;  i.e. it is an autonomous finite automaton.

``$\AFA\preceq \Turing_{\bL,\bR}$'': Conversely, if $\Ups:\sD\into\sD$
is an autonomous finite automaton, then let $\TM=(\sD,\tau,\Ups_*,\vV)$
be the $(\bL,\bR)$-Turing machine where
$\Ups:\sA\x\sD\x\sA\into\sA$ is defined by $\Ups_*(\zl,\zd,\zr):=\Ups(\zd)$,
and the functions $\vV$, $\tauL$, $\tauC$, and $\tauR$ are not important (for simplicity,
assume they are constants).
Then clearly, the `head dynamics' of 
$(\sD,\tau,\Ups_*,\vV)$ is an emulation of $\Ups:\sD\into\sD$.
\ethmprf

\Remarks{(a) If $h(\bL,\shift{})>0$  and $h(\bR,\shift{})>0$,
then Propositions \ref{Turing.CA}(b) and
 \ref{machine.equivalence}(a) imply that some questions about
the long-term behaviour of an $(\bL,\bR)$-defect particle are formally
undecidable.  For example, the question of whether the defect
particle eventually stops moving is equivalent to the Halting Problem.
Sutner \cite{Sutner} has identified similar undecidability 
issues for defect behaviour.

(b) The $(\betL,\betR)$-encoding mechanism in 
Proposition \ref{machine.equivalence} is quite crude;  a much more
efficient encoding could be obtained using finite state codes \cite[Ch.5]{LindMarcus}.

(c) In the standard definition, a Turing machine tape has
only a finite segment of nontrivial information;  we do
not assume this.  Likewise, in a standard pushdown automaton, the
stack has finite (but unbounded) height, whereas our
definition allows an infinitely high stack.

(d) Let $\TM$ be an APDA.  When moving
to the right (i.e. into the $\bR$-segment),  $\TM$ acts like a finite
automaton with state domain $\sD$, reading an $\sA$-valued
input stream and producing no output.  When moving to the left
(i.e. into the constant $\bar0$-segment), $\TM$ acts like an
autonomous finite automaton with state domain $\sD\x\sA$ and
update rule $\tl\Ups:\sD\x\sA\into\sD\x\sA$ 
defined by $\tl\Ups(\zd,\zr):=( \Ups(0,\zd,\zr),\tau(\zd,\zr))$.
A {\dfn runaway cycle} for $\TM$
is an $\tl\Ups$-periodic orbit $\{(\zd_p,\zr_p)\}_{p=1}^{P}$ [i.e. 
$\tl\Ups(\zd_p,\zr_p)=(\zd_{p+1},\zr_{p+1})$ and $\tl\Ups(\zd_P,\zr_P)=(\zd_{1},\zr_{1})$]
such that $\vV(\zd_p,\zr_p)=-1$ for all $p\in\CC{1...P}$.  In this case,
$\TM$ moves leftwards forever, and essentially belongs to the {\em Ballistic}
regime of \S\ref{S:ballistic}.  Not every APDA has a runaway cycle.
However, a variation of the Pumping Lemma
shows that, if $\TM$ moves leftward for long enough, it must
enter a runaway cycle.  Also, if $\TM$ has a runaway cycle which is
reachable from any initial conditions, and if $\rho\in\sM(\bR^+)$ is
a Bernoulli measure with full support, then for
 $\forall_\rho \ \br\in\bR^+$, an APDA with stack $\br$
will eventually enter a runaway state (see \S\ref{S:diffuse}
for definitions of `Bernoulli' and `$\forall_\rho$').

(e) By combining the arguments of  Theorem \ref{defect.random.walk}
and Propositions \ref{Turing.CA}(a) and \ref{machine.equivalence}(b), we can show that a defect
with a $\Phi$-fixed domain on one side and a $\Phi$-resolving subshift
on the other behaves like a pushdown automaton driven by a Markov process.
This is the `Markov Pushdown Automaton' 
regime in Table \ref{table:defect.regime}.

(f) The obvious multidimensional analogy of  Theorem \ref{Turing.CA}
involves a multidimensional Turing machine \cite[\S7.5]{HopcroftUllman}.
However, the problem of encoding a multidimensional bit array 
using a multidimensional subshift of finite type
(analogous to the $(\betL,\betR)$-encoding mechanism in 
Proposition \ref{machine.equivalence}) becomes much more complex.

(g) A completely different mechanism for universal computation has
been implemented using the (ballistic) defect dynamics of ECA\#110;
see \cite{Cook}, \cite{McIntosh2} or \cite[Chap.11]{Wolfram2}.
}

\subsection*{Conclusion}
  We have described the propagation of defects under the
action of cellular automata, but
many questions remain.  For example, we assumed that the defects
remain bounded in size, and act like `particles', as is the case in
well-known examples such as ECAs \#54, \#62, \#110, and \#184.  In
general, however, defects may grow over time like `blights' which
invade the whole lattice.  What are necessary/sufficient conditions
for the defect to remain bounded?  (In general, this is probably
formally undecidable; see \cite[Thm.3.2]{Sutner}.)

  Our theory is limited to one-dimensional subshifts of finite type.
This excludes some important cases (such as ECA \#18), where the
invariant subshift is sofic.  Can our theory be extended to sofic
shifts?  (Eloranta's `invariant subalphabet' approach covers some
sofic shifts by passing to a higher power presentation; see
\cite{ElNu,Elo93a,Elo93b}). 

  Even when $\bL$ and $\bR$ are subshifts of finite type, we only
understand defect dynamics in the polar opposite cases of `extreme
order' (i.e. $\bL$ and/or $\bR$ are $\Phi$-periodic) and `extreme
chaos' (i.e. $\bL$ and/or $\bR$ are $\Phi$-resolving, and endowed with
Parry measures).  We have been conspicuously silent about the so-called
`complicated' regime in Table \ref{table:defect.regime}.  In this
regime, pretty much anything can happen.  To see this, let $\sL$ and
$\sR$ be two disjoint finite alphabets, and let $\Phi_\sL:\LZ\into\LZ$
and $\Phi_\sR:\RZ\into\RZ$ be any two cellular automata
with local rules of radius 1.  Let $\sA:=
\sL\disj\sR$, and let $\Phi:\AZ\into\AZ$ be any radius-1
cellular automata such that $\Phi\restr{\LZ}=\Phi_\sL$ and
$\Phi\restr{\RZ}=\Phi_\sR$. Let $\ba:=[\bl \ \br]$ where
$\bl\in\sL^{\OO{-\oo...0}}$ and $\br\in\sR^{\CO{0...\oo}}$; then $\ba$
has a zero-width $(\bL,\bR)$-defect, where $\bL:=\LZ$ and $\bR:=\RZ$.
This defect must persist over time and can move either left or right
with unit speed.  If $\Phi$ has local rule
$\phi:\sA^{\CC{-1..1}}\into\sA$, then the defect's next move is determined by
the restriction of $\phi$ to the set $\sA^{\CC{-1..1}}\setminus
(\sL^{\CC{-1..1}}\disj\sR^{\CC{-1..1}})$.  However, the long-term
behaviour of the defect also depends on the dynamics of the CA
$(\LZ,\Phi_\sL)$ and $(\RZ,\Phi_\sR)$, which determine the `input
signals' which drive the defect.  Thus, the defect's behaviour is
potentially at least as complicated as the dynamics of any
one-dimensional CA, which could be very complicated indeed.

  However, perhaps if we control the topological dynamics of
$(\LZ,\Phi_\sL)$ and $(\RZ,\Phi_\sR)$, we can extend the
classification of Table
\ref{table:defect.regime}.  For example, perhaps we could weaken the
assumption of `$\Phi$-periodic' to `equicontinuous' in the {\em
Ballistic} and machine-emulating regimes, or perhaps we could replace
`right/left-resolving' with `positively expansive' in the {\em
Diffusive} regime.  Also, if $(\LZ,\Phi_\sL)$ and $(\RZ,\Phi_\sR)$
themselves manifest emergent defect dynamics, then perhaps we can
analyze the behaviour of the $(\bL,\bR)$-defect through its
interaction with these other defect particles (just as the Brownian
motion of a macromolecule is driven by a continual bombardment of
micromolecules).

 Finally, can a comparable theory of defect
particle kinematics be developed for subshifts of $\ZD$ for $D>1$?
Higher-dimensional shifts also admit infinitely extended defects
shaped like `curves' or `surfaces' 
\cite{Elo95,PivatoDefect2,PivatoDefect1}; what sort of motion do they
exhibit?  A general theory is probably hopeless: even interface curves
in a two-dimensional boolean CA exhibit a bewildering variety and
complexity of behaviour
\cite[\S3-\S6]{Gr98b}.  However, perhaps some
special cases are tractable.

\paragraph*{Acknowledgements:} 
  This paper was written during a research leave at Wesleyan
University, and was partially supported by the Van Vleck Fund.  I am
grateful to Ethan Coven, Adam Fieldsteel, Mike and Mieke Keane, and
Carol Wood for their generosity and hospitality. This research was
also supported by NSERC.  Finally, I would like to thank the two
referees for many perceptive and helpful suggestions about the notation
and exposition.

{\footnotesize
\bibliographystyle{alpha}
\bibliography{bibliography}
}

\end{document}